\documentclass[pdflatex,sn-mathphys-num]{sn-jnl}

\usepackage[utf8]{inputenc}
\usepackage{amsmath,amssymb,amsfonts}
\usepackage{amsthm}
\usepackage{graphicx}
\usepackage{xcolor}
\usepackage{tikz}
\usetikzlibrary{arrows.meta,positioning}
\usepackage{cleveref} 

\theoremstyle{thmstyleone}
\newtheorem{theorem}{Theorem}[section]
\newtheorem{proposition}[theorem]{Proposition}
\newtheorem{lemma}[theorem]{Lemma}
\newtheorem{corollary}[theorem]{Corollary}

\theoremstyle{thmstylethree}
\newtheorem{definition}[theorem]{Definition}
\newtheorem{example}[theorem]{Example}

\theoremstyle{thmstyletwo}
\newtheorem{remark}[theorem]{Remark}

\begin{document}

\title{The Isomorphism Problem for Cayley Graphs on Monogenic Semigroups}

\author[1]{\fnm{Lu} \sur{Li}}\email{lulimath@163.com}

\author*[1]{\fnm{Yong} \sur{Shao}}\email{yongshaomath@gm126.com}
\equalcont{These authors contributed equally to this work.}

\author[1]{\fnm{Lingli} \sur{Zeng}}\email{zengll929@nwu.edu.cn}
\equalcont{These authors contributed equally to this work.}
\maketitle
\vspace{0.5cm}
\noindent \textbf{Statements and Declarations}

\noindent \textbf{Competing Interests:} The authors declare that they have no competing interests.

\noindent \textbf{Funding:} Lu Li, Yong Shao, and Lingli Zeng are supported by Shaanxi provincial Natural Science Basic Research Program Project (No.2025JC-YBMS-086).

\noindent \textbf{Data Availability:} Data sharing is not applicable to this article as no datasets were generated or analysed during the current study.

\noindent \textbf{Acknowledgments:} The authors would like to thank Yong Shao Professor for his valuable suggestions on the early draft of this manuscript. 

\abstract{In this paper, we investigate the structural properties of Cayley graphs of monogenic semigroups and present a complete characterization of their weakly connected components. Several essential properties of these components are rigorously proven. Building upon these structural results, we introduce the notion of labelled colour-permuting isomorphisms between Cayley graphs. Using this new isomorphism framework, we address the Cayley graph isomorphism problem for monogenic semigroups, commonly known as the CI problem. In particular, we establish the necessary and sufficient conditions under which two Cayley graphs over monogenic semigroups are labelled colour-permuting isomorphic.}


\maketitle

\section{Introduction and preliminaries}

Cayley graphs play an important role in the study of groups and semigroups. Godsil and Royle~\cite{GodsilRoyle2001} established many important results concerning Cayley graphs of groups. Konstantinova and Lytkina~\cite{KonstantinovaLytkina2020} proved that Cayley graphs on finite groups with normal generating sets are integral under certain conditions. Pak and Radoicic~\cite{Pak2009} provided a detailed exposition of the Hamiltonicity problem for Cayley graphs of groups.  Kelarev and Praeger~\cite{KelarevPraeger2003} researched transitive Cayley graphs associated with semigroups. Kelarev~\cite{Kelarev2006} investigated Cayley graphs of inverse semigroups. Luo et al.~\cite{Luo2011} studied Cayley graphs of completely simple semigroups. Yang et al.~\cite{Yang2026} also provided a detailed discussion of the relationship between Cayley graphs and power graphs on completely 0-simple semigroups. 
\par The Cayley graph isomorphism problem is a fundamental topic in algebraic graph theory. Based on isomorphisms of Cayley graphs, Ádám~\cite{Adam1967} proposed the conjecture that all cyclic groups are CI‑groups. Babai~\cite{Babai1977} transformed the Cayley graph isomorphism problem to the problem of checking whether regular subgroups of the automorphism group are conjugate. Muzychuk~\cite{Muzychuk1995} proved that Ádám's conjecture holds in the square‑free order case in circulant graphs, i.e. Cayley graphs of cyclic groups and he~\cite{Muzychuk1997} obtained further related results. More generally, the isomorphism problem for Cayley graphs of other finite groups has been extensively investigated, for example, Li~\cite{Li1998} studied the isomorphism problem for connected Cayley graphs, Nowitz~\cite{Nowitz1992}, Xu~\cite{Xu1998}, and Hirasaka and Muzychuk~\cite{HirasakaMuzychuk2001} obtained both positive and negative results on the CI problem for abelian groups.
\par As the simplest finite semigroups generated by a single element, finite monogenic semigroups serve as a vital bridge between cyclic groups and general finite semigroups. Extensive studies on monogenic semigroups have laid a fundamental framework for the classification and enumeration of finite commutative semigroups, and yielded abundant elementary examples for applications in automata theory, combinatorics, and cryptography. Das et al.~\cite{Das2013} studied an undirected graph associated with a finite monogenic semigroup with zero, and computed its graph parameters, including diameter, clique number, and chromatic number. 
\par In this paper, we study Cayley graphs of finite monogenic semigroups and investigate their isomorphism problem through this index-period structure. We first describe the weakly connected components of such Cayley graphs in terms of the arithmetic data of the connection set, especially the greatest common divisor of its exponents with the period. We then introduce labelled colour-permuting isomorphisms between Cayley graphs. This notion generalizes colour-preserving automorphisms and gives a suitable framework for formulating the CI problem in the monogenic semigroup setting.

Our main results give necessary and sufficient conditions for two Cayley graphs over a finite monogenic semigroup to be labelled colour-permuting isomorphic. The analysis is divided according to the location of the connection sets: connection sets contained in the non-cyclic part, connection sets contained in the cyclic part, and mixed connection sets. In the cyclic case, the criteria are expressed through admissible ordered-reduced connection sets, gcd conditions, and rotation data on the corresponding cycles. In the mixed case, the non-cyclic part is shown to impose rigid constraints, while the cyclic part controls the remaining isomorphism behaviour.

The paper is organized as follows. In Section 2, we define labelled colour-permuting isomorphisms and formulate the CI problem for finite monogenic semigroups. In Section 3, we determine the weakly connected components of Cayley graphs of finite monogenic semigroups and establish their basic structural properties. In Section 4, we use these structural results to classify labelled colour-permuting isomorphisms for the three types of connection sets described above.

A finite monogenic semigroup generated by $a$ has the following representation:
\begin{align*}
    \left\langle a \right\rangle=\langle a\mid a^{m+r}=a^m\rangle
    =\{a,a^2,\ldots, a^{m-2}, a^{m-1},a^m,\ldots,a^{m+r-1}\},
\end{align*}
where $m$ is the \emph{index} and $r$ is the \emph{period}. Thus, for every $s\ge m$,
\begin{align*}
    a^{s+r}=a^s.
\end{align*}
\par We divide $\left\langle a \right\rangle$ into two subsets, 
\begin{align*}
    \{a,a^2,\ldots, a^{m-2}, a^{m-1}\} \text{ and }
    \{a^m,a^{m+1},\ldots,a^{m+r-1}\}.
\end{align*}
We call the elements of $\{a,a^2,\ldots, a^{m-2}, a^{m-1}\}$ non‑cyclic elements and this set is called the non‑cyclic part of the monogenic semigroup, denoted by $\textbf{\textit{NCP}}$, the elements of $\{a^m,a^{m+1},\ldots,a^{m+r-1}\}$ cyclic elements and this set is called the cyclic part of the monogenic semigroup, denoted by $\textbf{\textit{CP}}$. The cyclic part behaves similarly to a cyclic group, while the non-cyclic part records genuinely semigroup-theoretic behaviour. Therefore, finite monogenic semigroups provide a natural testing ground for extending the CI problem from cyclic groups to semigroups.
\par Unless otherwise specified, all monogenic semigroups in the remainder of this paper are finite and are written in the form $\left\langle a \right\rangle$ with index $m$ and period $r$.

Let $S$ be a semigroup and let $C$ be a subset of $S$. The Cayley graph $\text{Cay}(S, C)$ is the graph with vertex set $S$ and arc set
\begin{align*}
    E\bigl(\text{Cay}(S,C)\bigr)
    =\{(s,cs):s\in S,\ c\in C\}.
\end{align*}
For arbitrary $a^u,a^v \in \left\langle a \right\rangle$, $a^u\to a^v$ represents $(a^u, a^v) \in E\bigl(\text{Cay}(\left\langle a \right\rangle,C)\bigr)$; $a^u \gets a^v$ represents $(a^v, a^u) \in E\bigl(\text{Cay}(\left\langle a \right\rangle,C)\bigr)$.
\par Two directed graphs $\Gamma_1$ and $\Gamma_2$ are isomorphic, written $\Gamma_1 \cong \Gamma_2$, if there is a bijection $f:V(\Gamma_1)\to V(\Gamma_2)$ such that
\begin{align*}
    (u,v)\in E(\Gamma_1)
    \Longleftrightarrow
    (f(u),f(v))\in E(\Gamma_2).
\end{align*}
\par Let both $H$ and $C$ be subsets of $\left\langle a \right\rangle$, the induced subgraph of $H$ in $\text{Cay}(\left\langle a \right\rangle, C)$ is denoted by $\text{Cay}(\left\langle a \right\rangle, C)[H]$.

We shall repeatedly use the elementary fact that for any positive integer $i, r$ the congruence equation $ix\equiv0\pmod r$ has a least positive solution $\frac{r}{\gcd(i,r)}$ \cite{Shoup2009}. Based on this fact, the following simple conclusion can be obtained. For arbitrary connection set $C\subseteq \left\langle a \right\rangle$ and $a^i \in C$, let $a^u \in \textbf{\textit{CP}}$,
\begin{align*}
    a^u \to a^{u+i} \to \cdots \to a^{u+(\frac{r}{\gcd(i,r)}-1)i} \to a^u
\end{align*}
is a directed circle in $\text{Cay}(\left\langle a \right\rangle, C)$. 

\section{CI problem on a finite monogenic semigroup}

In this section, we set up the monogenic-semigroup version of the CI problem. We first recall the index-period decomposition and the non-cyclic/cyclic split, then introduce labelled colour-permuting isomorphisms and the corresponding notion of CI-object. These notions provide the language needed for the structural results in the next section.

\begin{proposition}\label{pro: c to c}
    If there exist two nonempty subsets $C_1 ,C_2 \subseteq \left\langle a \right\rangle$, such that $\text{Cay}(\left\langle a \right\rangle,C_1) \cong \text{Cay}(\left\langle a \right\rangle,C_2)$, then for any graph isomorphism $f$ from $\text{Cay}(\left\langle a \right\rangle,C_1)$ to $\text{Cay}(\left\langle a \right\rangle,C_2)$, we have that  
\begin{gather*}
    (\forall a^u \in \textbf{\textit{NCP}}) \text{ } f(a^u) \in \textbf{\textit{NCP}}, \\
    (\forall a^v \in \textbf{\textit{CP}}) \text{ } f(a^v) \in \textbf{\textit{CP}}.
\end{gather*}
\end{proposition}

\begin{proof}
    Suppose that $f$ is a graph isomorphism from $\text{Cay}(\left\langle a \right\rangle,C_1)$ to $\text{Cay}(\left\langle a \right\rangle,C_2)$ and $a^u, a^v \in \left\langle a \right\rangle$ such that $f(a^u) = a^v$. If $a^u \in \textbf{\textit{CP}}$, but $a^v \in \textbf{\textit{NCP}}$, then for any $a^i \in C_1$, there exists a directed circle in $\text{Cay}(\left\langle a \right\rangle, C_1)$,
    \begin{align*}
        a^u \to a^{u+i} \to \cdots \to a^{u+(\frac{r}{\gcd(i,r)}-1)i} \to a^u.
    \end{align*}
    Correspondingly, $\text{Cay}(\left\langle a \right\rangle, C_2)$ also contains directed circle,
    \begin{align*}
        f(a^u) \to f(a^{u+i}) \to \cdots \to f(a^{u+(\frac{r}{\gcd(i,r)}-1)i}) \to f(a^u),
    \end{align*}
    where $f(a^u) = a^v$ and $a^v$ is a non-cyclic element. Since no \textbf{\textit{NCP}} element lies on a directed cycle, it follows that $a^v \in \textbf{\textit{CP}}$. Therefore, it is clear that if $a^u \in \textbf{\textit{NCP}}$, then $a^v \in \textbf{\textit{NCP}}$. 
\end{proof}

\begin{proposition}\label{pro: equivant}
Suppose that $C_1$ and $C_2$ are distinct singleton subsets of $\left\langle a \right\rangle$. Then $\text{Cay}(\left\langle a \right\rangle, C_1) = \text{Cay}(\left\langle a \right\rangle, C_2)$ if and only if $C_1 \cup C_2 = \{a^{m-1}, a^{m+r-1}\}$.
\end{proposition}
\begin{proof}
For any $a^u \in \left\langle a \right\rangle$, we have $a^u \cdot a^{m-1} = a^u \cdot a^{m+r-1}$. It is clear that $\text{Cay}(\left\langle a \right\rangle, \{a^{m-1}\}) = \text{Cay}(\left\langle a \right\rangle, \{a^{m+r-1}\})$.
\par Conversely, suppose that there are $\left\{a^i\right\} \subseteq \left\langle a \right\rangle$ and $\left\{a^j\right\} \subseteq \left\langle a \right\rangle$ such that $\text{Cay}(\left\langle a \right\rangle, \{a^i\}) = \text{Cay}(\left\langle a \right\rangle, \{a^j\})$. Then $a^i \cdot a = a^j \cdot a$. Since $C_1 \ne C_2$, $i \ne j$, $i = j + kr$ or $j = i + kr$, and $i, j \ge m-1$. Hence, $i = m-1, j = m + r - 1$ or $i = m+r-1, j = m-1$. Thus $C_1 \cup C_2 = \{a^{m-1}, a^{m+r-1}\}$.
\end{proof}

\begin{definition}\label{def: rc set}
    Define an operator $\Psi$ on the power set $2^{ \left\langle a \right\rangle }$ as follows:
\begin{align*}
    \Psi(C)=
    \begin{cases}
    (C\setminus\{a^{m-1}\})\cup\{a^{m+r-1}\},&a^{m-1}\in C,\\
    C,&a^{m-1}\notin C.
    \end{cases}
\end{align*}
And for a subset $C \subseteq \left\langle a \right\rangle$, $\Psi(C)$ is called the $reduced \text{ } connection \text{ } set$ of $C$.
\end{definition}

By Proposition~\ref{pro: equivant}, 
\begin{align*}
    \text{Cay}(\left\langle a \right\rangle,C)=\text{Cay}(\left\langle a \right\rangle, \Psi(C)).
\end{align*}
Unless stated otherwise, connection sets in the remainder of the paper are reduced.

\begin{proposition}\label{pro: connection set}
 If there exist two nonempty subsets $C_1, C_2 \subseteq \left\langle a \right\rangle$, such that $\text{Cay}(\left\langle a \right\rangle,C_1) \cong \text{Cay}(\left\langle a \right\rangle,C_2)$, then $|C_1| = |C_2|$.
\end{proposition}
\begin{proof}
    Let $d^+(x)$ denote the out-degree of $x$. Let $\Delta^+(\text{Cay}(\left\langle a \right\rangle,C))$ denote the maximum out-degree of all vertices in $\text{Cay}(\left\langle a \right\rangle,C)$. It is clear that $d^+(x) \le |C_1|$ for any $x \in \left\langle a \right\rangle$ and $\Delta^+(\text{Cay}(\left\langle a \right\rangle,C_1)) \le |C_1|$. Suppose that $a^p,a^q \in C_1 $ and $p \ne q$ satisfying $a^p \cdot a = a^q \cdot a$. Since $p \ne q$, $p \ge m-1 , q \ge m-1$ and $p \equiv q \pmod r$. Thus $p = m-1, q = m+r-1$ or $p = m+r-1, q = m-1$. But $a^{m-1}$ is neither in $C_1$ nor in $C_2$. Thus $d^+(a) = |C_1| \text{ and } \Delta^+(\text{Cay}(\left\langle a \right\rangle,C_1)) = |C_1|$. Similarly, $\Delta^+(\text{Cay}(\left\langle a \right\rangle,C_2)) = |C_2|$, and so, $|C_1| = |C_2|$.
\end{proof}

\begin{definition}\label{def: color-preserving automorphism}
(\cite{Khosravi2021})
Let $\left\langle a \right\rangle$ be a monogenic semigroup and $C$ be a subset of $\left\langle a \right\rangle$. A bijection $f$ is called a colour-preserving automorphism if $xc = y$ implies $f(x)c = f(y)$ for every $x,y \in \left\langle a \right\rangle$ and $c \in C$.
\end{definition}

\begin{definition}\label{def: labelled color-permuting isomorphic}
Let $\left\langle a \right\rangle$ be a monogenic semigroup. $C_1, C_2 \subseteq \left\langle a \right\rangle$, $C_1,C_2 \ne \emptyset $. $\text{Cay}(\left\langle a \right\rangle,C_1) $ and $ \text{Cay}(\left\langle a \right\rangle,C_2)$ are \textit{labelled colour-permuting isomorphic}, if there exist a bijection $f : \left\langle a \right\rangle \to \left\langle a \right\rangle$ and a bijection $\varphi : C_1 \to C_2$, such that $f(xa) = \varphi(x)f(a)$ for all $a \in \left\langle a \right\rangle$ and all $x \in C_1$.
\par We write the labelled colour‑permuting isomorphism between $\text{Cay}(\left\langle a \right\rangle,C_1)$ and $\text{Cay}(\left\langle a \right\rangle,C_2)$ as $\text{Cay}(\left\langle a \right\rangle,C_1) \simeq \text{Cay}(\left\langle a \right\rangle,C_2)$. If $\text{Cay}(\left\langle a \right\rangle,C_1) \simeq \text{Cay}(\left\langle a \right\rangle, \allowbreak C_2)$, it is easy to see that $\text{Cay}(\left\langle a \right\rangle,C_1) \cong \text{Cay}(\left\langle a \right\rangle,C_2) $. 
\end{definition}

By Definitions~\ref{def: color-preserving automorphism} and~\ref{def: labelled color-permuting isomorphic}, the colour‑preserving automorphism in Definition~\ref{def: color-preserving automorphism} is a special case of Definition~\ref{def: labelled color-permuting isomorphic}.

In the following, we shall use the labelled colour‑permuting isomorphism to define the CI problem for monogenic semigroups.

\begin{definition}\label{def: CI}
Let $\left\langle a \right\rangle$ be a monogenic semigroup. The Cayley graph $\text{Cay}(\left\langle a \right\rangle,\allowbreak C)$ is called a \textit{CI-object} of $\left\langle a \right\rangle$, if all Cayley graphs of $\left\langle a \right\rangle$ that are isomorphic to $\text{Cay}(\left\langle a \right\rangle, C)$ are also labelled colour-permuting  isomorphic to $\text{Cay}(\left\langle a \right\rangle, C)$. $\left\langle a \right\rangle$ is called a \textit{CI monogenic semigroup} if all Cayley graphs of $\left\langle a \right\rangle$ are CI-objects.
\end{definition}

It should be noted that not all Cayley graphs on $\left\langle a \right\rangle$ are CI‑objects. To illustrate this fact, we shall present the following example.

\begin{example}
Let $\left\langle a \right\rangle$ be a monogenic semigroup with index 2 and periodic 8. Take $C_1 = \left\{a^3, a^6, a^7\right\}$, $C_2 = \left\{a^2, a^3, a^7\right\}$. Then $\text{Cay}(\left\langle a \right\rangle,C_1) \cong \text{Cay}(\left\langle a \right\rangle,C_2)$ but $\text{Cay}(\left\langle a \right\rangle,C_1) \not\simeq \text{Cay}(\left\langle a \right\rangle,C_2)$. Thus $\text{Cay}(\left\langle a \right\rangle,C_1)$ and $\text{Cay}(\left\langle a \right\rangle,C_2)$ are not CI-object.
\end{example}
\begin{proof}
Define a bijection $f: \left\langle a \right\rangle \to \left\langle a \right\rangle$ by
\begin{align*}
    \begin{array}{c|ccccccccccc}
        x&a&a^2&a^3&a^4&a^5&a^6&a^7&a^8&a^9\\  \hline
        f(x)&a&a^2&a^7&a^8&a^5&a^6&a^3&a^4&a^9
    \end{array}
\end{align*}
    These two Cayley graphs are of the following form
    
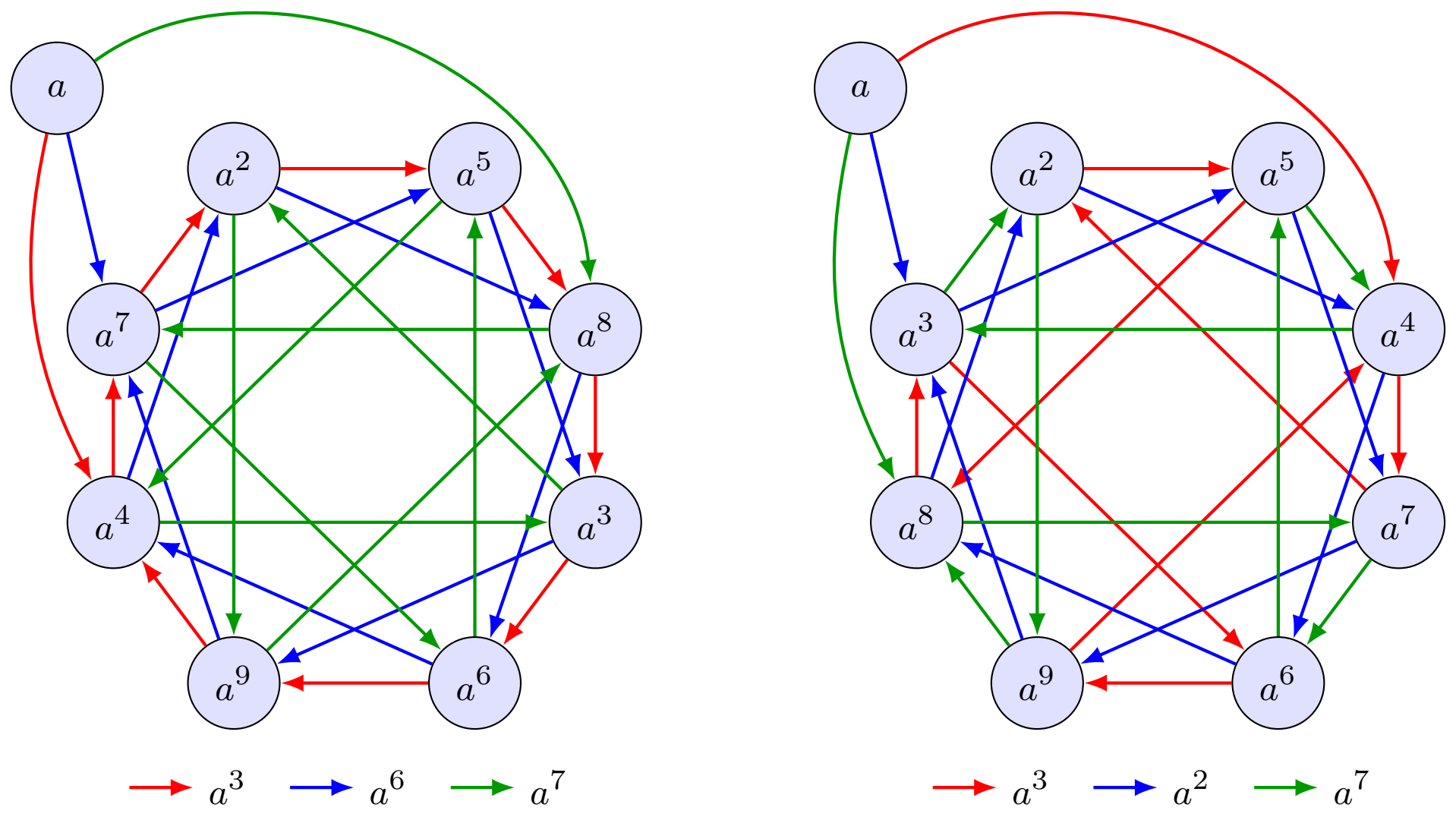
\begin{figure}[h]
    \centering
\begin{tikzpicture}[scale=0.7,
     vertex/.style={
        circle,
        draw,
        fill=blue!12,
        minimum size=8mm,
        inner sep=0pt
    },
    gen3/.style={-{Latex[length=2mm]}, thick, red},
    gen2/.style={-{Latex[length=2mm]}, thick, blue},
    gen6/.style={-{Latex[length=2mm]}, thick, blue},
    gen7/.style={-{Latex[length=2mm]}, thick, green!60!black},
     every node/.style={font=\small}
]
\begin{scope}[shift={(-8,0)}, local bounding box=first]
    \node[vertex] (v1)  at (-3.2,4.2)  {$a$};
    \node[vertex] (v2)  at (-1,3.2)   {$a^2$};
    \node[vertex] (v3)  at (3.5,-1.2)   {$a^3$};
    \node[vertex] (v4)  at (-2.5,-1.2)  {$a^4$};
    \node[vertex] (v5)  at (2,3.2)   {$a^5$};
    \node[vertex] (v6)  at (2,-3.2)   {$a^6$};
    \node[vertex] (v7)  at (-2.5,1.2)   {$a^7$};
    \node[vertex] (v8)  at (3.5,1.2)   {$a^8$};
    \node[vertex] (v9)  at (-1,-3.2)   {$a^9$};
            
    \draw[gen3] (v1)  to[bend right=20] (v4);
    \draw[gen3] (v2)  -- (v5);
    \draw[gen3] (v5)  -- (v8);
    \draw[gen3] (v8)  -- (v3);
    \draw[gen3] (v3)  -- (v6);
    \draw[gen3] (v6)  -- (v9);
    \draw[gen3] (v9)  -- (v4);
    \draw[gen3] (v4)  -- (v7);
    \draw[gen3] (v7)  -- (v2);
    
    \draw[gen6] (v1)  -- (v7);
    \draw[gen6] (v2)  -- (v8);
    \draw[gen6] (v8)  -- (v6);
    \draw[gen6] (v6)  -- (v4);
    \draw[gen6] (v4)  -- (v2);
    \draw[gen6] (v5)  -- (v3);
    \draw[gen6] (v3)  -- (v9);
    \draw[gen6] (v9)  -- (v7);
    \draw[gen6] (v7)  -- (v5);
    
    \draw[gen7] (v1)  to[bend left=60] (v8);
    \draw[gen7] (v2)  -- (v9);
    \draw[gen7] (v9)  -- (v8);
    \draw[gen7] (v8)  -- (v7);
    \draw[gen7] (v7)  -- (v6);
    \draw[gen7] (v6)  -- (v5);
    \draw[gen7] (v5)  -- (v4);
    \draw[gen7] (v4)  -- (v3);
    \draw[gen7] (v3)  -- (v2);
        
    \begin{scope}[shift={(-2.3,-4.5)}]
        \draw[gen3] (0,0) -- +(0.8,0) node[right,black] {$a^3$};
        \draw[gen6] (2,0) -- +(0.8,0) node[right,black] {$a^6$};
        \draw[gen7] (4,0) -- +(0.8,0) node[right,black] {$a^7$};
    \end{scope}
\end{scope}

\begin{scope}[shift={(2,0)}, local bounding box=first]
\node[vertex] (v1)  at (-3.2,4.2)  {$a$};
    \node[vertex] (v2)  at (-1,3.2)   {$a^2$};
    \node[vertex] (v7)  at (3.5,-1.2)   {$a^7$};
    \node[vertex] (v8)  at (-2.5,-1.2)  {$a^8$};
    \node[vertex] (v5)  at (2,3.2)   {$a^5$};
    \node[vertex] (v6)  at (2,-3.2)   {$a^6$};
    \node[vertex] (v3)  at (-2.5,1.2)   {$a^3$};
    \node[vertex] (v4)  at (3.5,1.2)   {$a^4$};
    \node[vertex] (v9)  at (-1,-3.2)   {$a^9$};
            
    \draw[gen3] (v1)  to[bend left=60] (v4);
    \draw[gen3] (v2)  -- (v5);
    \draw[gen3] (v5)  -- (v8);
    \draw[gen3] (v8)  -- (v3);
    \draw[gen3] (v3)  -- (v6);
    \draw[gen3] (v6)  -- (v9);
    \draw[gen3] (v9)  -- (v4);
    \draw[gen3] (v4)  -- (v7);
    \draw[gen3] (v7)  -- (v2);
    
    \draw[gen2] (v1)  -- (v3);
    \draw[gen2] (v2)  -- (v4);
    \draw[gen2] (v4)  -- (v6);
    \draw[gen2] (v6)  -- (v8);
    \draw[gen2] (v8)  -- (v2);
    \draw[gen2] (v5)  -- (v7);
    \draw[gen2] (v7)  -- (v9);
    \draw[gen2] (v9)  -- (v3);
    \draw[gen2] (v3)  -- (v5);
    
    \draw[gen7] (v1)  to[bend right=20] (v8);
    \draw[gen7] (v2)  -- (v9);
    \draw[gen7] (v9)  -- (v8);
    \draw[gen7] (v8)  -- (v7);
    \draw[gen7] (v7)  -- (v6);
    \draw[gen7] (v6)  -- (v5);
    \draw[gen7] (v5)  -- (v4);
    \draw[gen7] (v4)  -- (v3);
    \draw[gen7] (v3)  -- (v2);
        
\begin{scope}[shift={(-2.3,-4.5)}]
    \draw[gen3] (0,0) -- +(0.8,0) node[right,black] {$a^3$};
    \draw[gen2] (2,0) -- +(0.8,0) node[right,black] {$a^2$};
    \draw[gen7] (4,0) -- +(0.8,0) node[right,black] {$a^7$};
\end{scope}
\end{scope}

\end{tikzpicture}%
\caption{The Cayley graph \(\operatorname{Cay}(\langle a\rangle, C_1)\) and \(\operatorname{Cay}(\langle a\rangle, C_2)\) in Example 2.8.}
    \label{fig:example2-8-C1}
\end{figure}

By figures 3 and 4, $f$ is a graph isomorphism.
\par Suppose that there exist a labelled colour-permuting isomorphism $g$ and a corresponding bijection $\varphi$. By Proposition~\ref{pro: c to c}, the unique non-cyclic vertex $a$ must be fixed by $g$. On the \textbf{\textit{CP}}, the colour $a^6$ produces directed cycles of length $4$, whereas $a^3$ and $a^7$ produce cycles of length $8$. Hence $\varphi(a^6)=a^2$, and $\varphi$ maps $\{a^3,a^7\}$ to itself.

If $\varphi(a^3)=a^3$, then $g(a^4)=a^4$ and $g(a^7)=a^3$. Applying the same identity with colour $a^3$ at $x=a^4$ yields
\begin{align*}
    g(a^7)=a^3g(a^4)=a^7.
\end{align*}
That is a contradiction. 

If $\varphi(a^3)=a^7$, then $g(a^4)=a^8$ and $g(a^7)=a^3$, now
\begin{align*}
    g(a^7)=a^7g(a^4)=a^7a^8=a^7.
\end{align*}
That is a contradiction. Therefore no labelled colour-permuting isomorphism exists.
\end{proof} 

\begin{definition}\label{def: cyclic CI}
(\cite{BabaiFrankl1979II})
   Let $G$ be a cyclic group. The Cayley graph $\text{Cay}(G, C)$ is called a CI-object of $G$, if all Cayley graphs of $G$ that are isomorphic to $\text{Cay}(G, C)$ are $\text{Cay}(G, \alpha(C))$ for some automorphism $\alpha$ on $G$. $G$ is called a CI group, if all Cayley graphs of $G$ are CI-objects. 
\end{definition}

Moreover, if the finite monogenic semigroup in Definition \ref{def: CI} is a cyclic group, then the equivalence of Definition \ref{def: CI} and Definition \ref{def: cyclic CI} can be established by the following theorem.

\begin{lemma}\label{lemma: Automorphism extended}
Let $G$ be a cyclic group and $H$ be a subgroup of $G$. For an arbitrary automorphism $\alpha$ of $H$, $\alpha$ can be extended to an automorphism on $G$.
\end{lemma}
\begin{proof}
    Let $G = \left\{g, g^2, \cdots, g^n \right\}$, $|G|= n$, $|H| = d$. Then $d|n$. If $d = 1$, it is clear that this conclusion is true. Assume that $d > 1$. Without loss of generality, let $H = \left\{g^{\frac{n}{d}}, g^{\frac{2n}{d}}, \cdots, g^n\right\}$. Suppose that $\alpha$ is an automorphism on $H$, and the value of $\alpha$ is uniquely determined by its value on the generator $g^{\frac{n}{d}}$. Let $\alpha(g^{\frac{n}{d}}) = g^{s\frac{n}{d}}$ and $\gcd(s, d) = 1$.
    \par Let $n = p_1^{e_1} p_2^{e_2} \cdots p_t^{e_t}$, where $p_1, p_2, \cdots, p_t$ are prime numbers. If $p_1|d$, then a congruence equation $b \equiv s \pmod {p_i^{e_i}}$ is obtained. If $p_1 \nmid d$, then a congruence equation $b \equiv 1 \pmod {p_i^{e_i}}$ is obtained. In this way, we obtain a system of congruence equations:
    \begin{align*}
        \begin{cases}
        b \equiv s \pmod {p_1^{e_1}}\\
        b \equiv 1 \pmod {p_2^{e_2}}\\
        \quad \quad \quad \cdots \\
        b \equiv s \pmod {p_t^{e_t}}
    \end{cases}
    \end{align*}
    Since $p_1^{e_1}, p_2^{e_2}, \cdots, p_t^{e_t}$ are coprime. By the $Chinese \text{ } remainder \text{ } theorem$, the solution $b$ of the congruence equation system can be obtained.
    \par If $b \equiv s \pmod {p_u^{e_u}}$ or $b \equiv 1 \pmod {p_u^{e_u}}$, then it can be obtained that $b \equiv s \pmod {p_u}$ or $b \equiv 1 \pmod {p_u}$ for all $u \in \left\{1, 2, \cdots, t\right\}$. Thus $\gcd(n, b) = 1$. Let $d = p_{i_1}^{e_{i_1}} p_{i_2}^{e_{i_2}} \cdots p_{i_t}^{e_{i_t}}$. Then $b \equiv s \pmod {p_{i_1}^{e_{i_1}}}, b \equiv s \pmod {p_{i_2}^{e_{i_2}}}, \cdots, b \equiv s \pmod {p_{i_t}^{e_{i_t}}}$. Thus $b \equiv s \pmod d$.
    \par Define $\psi :G\to G$, $\psi(g^k) = g^{bk}$. Since $\gcd(b, n) = 1$, $\psi$ is an automorphism on $G$. $\psi(g^{\frac{n}{d}\cdot t}) = g^{\frac{n}{d}\cdot t\cdot b} = g^{\frac{n}{d}\cdot t\cdot s} = \alpha(g^{\frac{n}{d}\cdot t})$. Thus $\psi|_H = \alpha$, $\psi$ is an extension of $\alpha$ to $G$. 
\end{proof}

\begin{theorem}\label{throrem: 1}
Let $G$ be a cyclic group and $C_1, C_2 \subseteq G$. Then $\text{Cay}(G,C_1) \simeq \text{Cay}(G, C_2)$ if and only if there exists an automorphism $\alpha$ of $G$ such that $\alpha(C_1) = C_2$.
\end{theorem}
\begin{proof} 
    Suppose that $f$ is the labelled colour-permuting isomorphism and $\varphi$ is the bijection from $C_1$ to $C_2$. Let $e$ is the identity element of $G$, and let $c = f(e)$. 
    \par Define $\sigma : G \to G$, $\sigma(a) = c^{-1}f(a)$. For any $a^u \in G$, $c\cdot a^u \in G$. Since $f$ is a bijection, there exists $a^v \in G$ such that $f(a^v) = c\cdot a^u \in G$, $\sigma(a^v) = a^u$. Thus $\sigma$ is a surjection. Suppose that $a^u, a^v \in G$ and $a^u \ne a^v$. Then $f(a^u) \ne f(a^v)$, and $c^{-1}f(a^u) \ne c^{-1}f(a^v)$. Thus $\sigma(a^u) \ne \sigma(a^v)$ and $\sigma$ is an injection. It is clear that $\sigma$ is a bijection, and $\sigma(e) = e$. For all $x \in C_1, a \in G$,
    \begin{align*}
        \sigma(xa) = c^{-1}f(xa) = c^{-1}\varphi(x)f(a) = c^{-1}\varphi(x)cc^{-1}f(a) = c^{-1}\varphi(x)c\sigma(a)
    \end{align*}
    Thus $\sigma(xa) = c^{-1}\varphi(x)c\sigma(a) = \sigma(x)\sigma(a)$. For any $x \in C_1$, $\sigma(x) = c^{-1}f(x) = c^{-1}f(x \cdot e) = c^{-1}\varphi(x)f(e) = c^{-1}\varphi(x)c = \varphi(x)$. Thus $\sigma(C_1) = C_2$. It is clear that $\sigma(\left\langle C_1 \right\rangle) = \left\langle C_2 \right\rangle$ and $|\left\langle C_1 \right\rangle| = |\left\langle C_2 \right\rangle|$, hence $\left\langle C_1 \right\rangle = \left\langle C_2 \right\rangle$.
    This shows that $\sigma|_{\left\langle C_1 \right\rangle}$ is an automorphism of $\left\langle C_1 \right\rangle$. By Lemma~\ref{lemma: Automorphism extended}, $\sigma|_{\left\langle C_1 \right\rangle}$ can extend to an automorphism $\alpha$ of $G$. And $\alpha(C_1) = C_2$. 
\end{proof}

\section{The structure of Cayley graphs}

In this section, we now study the structure of Cayley graphs of finite monogenic semigroups through the index-period decomposition. The non-cyclic part and the cyclic part behave differently, and the arithmetic data of the connection set, especially its gcd with the period, controls the weakly connected components of the graph. This description will later allow us to compare different components and derive the isomorphism criteria in Section 4.
In the following discussion, for a subset $C = \{a^{i_1}, a^{i_2}, \cdots, a^{i_n}\}$, let $\gcd(C, r)$ denote the greatest common divisor of the exponents of all elements in $C$ and the period $r$.

\begin{definition}\label{def: rho1 rho2}
Let $\left\langle a \right\rangle$ be a monogenic semigroup.
$C = \left\{a^{i_1}, a^{i_2}, \cdots, a^{i_n}\right\} \subseteq \left\langle a \right\rangle$, $d = \gcd(C, r)$. Define respecitively the relation $\rho_1$ on $\textbf{\textit{CP}}$ and $\rho_2$ on $\textbf{\textit{NCP}}$ as folows: 
\begin{gather*}
    \rho_1 = \{(a^{k_1}, a^{k_2}) \in \textbf{\textit{CP}} \times \textbf{\textit{CP}} | k_1 \equiv k_2\pmod{d}\}, \\ 
    \rho_2 = \{(a^{k_1}, a^{k_2}) \in \textbf{\textit{NCP}} \times \textbf{\textit{NCP}} | \\(\exists a^p, a^q \in C \text{ and } t_1,t_2 \in \mathbb{N}), (a^{t_1p+ k_1}, a^{t_2q+ k_2}) \in \rho_1 \}.
\end{gather*}
\end{definition}

\begin{proposition}\label{pro: equivalent relation}
Both $\rho_1$ and $\rho_2$ are equivalence relations.
\end{proposition}
\begin{proof}
    For any $a^{k_1} \in \textbf{\textit{CP}}$, since $k_1 \equiv k_1 \pmod d$, $(a^{k_1},a^{k_1}) \in \rho_1$.
    For any $a^{k_1}, a^{k_2} \in \textbf{\textit{CP}}$, if $(a^{k_1}, a^{k_2}) \in \rho_1$, then $k_1 \equiv k_2 \pmod d$, so $k_2 \equiv k_1 \pmod d$. Thus $(a^{k_2}, a^{k_1}) \in \rho_1$.
    For any $a^{k_1}, a^{k_2}, a^{k_3} \in \textbf{\textit{CP}}$, if $(a^{k_1}, a^{k_2}) \in \rho_1$ and $(a^{k_2}, a^{k_3}) \in \rho_1$, then $k_1 \equiv k_2 \pmod d$ and $k_2 \equiv k_3 \pmod d$ and so $k_1 \equiv k_3 \pmod d$. Thus $(a^{k_1}, a^{k_3}) \in \rho_1$ and so $\rho_1$ is an equivalent relation. 
    \par For any $a^{k_1} \in \textbf{\textit{NCP}}$, there exist $t_1 \in \mathbb{N}$ and $a^p \in C$, such that $a^{t_1p+k_1} \in \textbf{\textit{CP}}$ and $t_1p+k_1 \equiv t_1p+k_1 \pmod d$, so $(a^{k_1}, a^{k_1}) \in \rho_2$.
    For any $a^{k_1}, a^{k_2} \in \textbf{\textit{NCP}}$, if $(a^{k_1}, a^{k_2}) \in \rho_2$, then $(a^{t_1p+k_1}, \allowbreak a^{t_2q+k_2}) \in \rho_1$ for some $a^p, a^q \in C$ and $t_1, t_2 \in \mathbb{N}$, so $(a^{t_2q+k_2}, a^{t_1p+k_1}) \in \rho_1$. Thus $(a^{k_2}, a^{k_1}) \in \rho_2$.
    For arbitrary $a^{k_1}, a^{k_2}, a^{k_3} \in \textbf{\textit{NCP}}$, if $(a^{k_1}, a^{k_2}) \in \rho_2$ and $(a^{k_2}, a^{k_3}) \in \rho_2$, then $(a^{t_1p_1+k_1}, a^{t_2q_1+k_2}) \in \rho_1 $ and $ (a^{t_3p_2+k_2}, a^{t_4q_2+k_3}) \in \rho_1$ for some $a^{p_1}, a^{p_2}, a^{q_1}, a^{q_2} \in C$ and $t_1, t_2, t_3, t_4 \in \mathbb{N}$. So $t_1p_1+k_1 \equiv t_2q_1+k_2 \pmod d, t_3p_2+k_2 \equiv t_4q_2+k_3 \pmod d$. Since $d = \gcd(i_1, i_2, \cdots, i_n, r)$ and $a^{p_2}, a^{q_1} \in C$, then $t_3p_2 \equiv t_2q_1 \pmod d$. So $t_2q_1+k_2 \equiv t_3p_2+k_2 \pmod d$, then $t_1p_1+k_1 \equiv t_4q_2+k_3 \pmod d$. Thus $(a^{k_1}, a^{k_3}) \in \rho_2$ and so $\rho_2$ is an equivalent relation.
\end{proof}

\begin{proposition}\label{pro: weakly connected}
For arbitrary $a^{k_1}, a^{k_2} \in \textbf{\textit{CP}}$, $(a^{k_1}, a^{k_2}) \in \rho_1$ if and only if $a^{k_1}$ and $ a^{k_2}$ are weakly connected in $\text{Cay}(\left\langle a \right\rangle, C)$.
\end{proposition}
\begin{proof}
    Suppose that $a^{k_1}, a^{k_2} \in \textbf{\textit{CP}}$ and $(a^{k_1}, a^{k_2}) \in \rho_1$. Then $k_1 \equiv k_2 \pmod d$, $d|(k_1-k_2)$. Let $k_2-k_1 = sd$ for some $s >0$. Since the $\textit{B\'ezout } Lemma$, there exist $x_1, x_2, \cdots, x_n, \bar{x} \in \mathbb{Z}$, such that $d = x_1i_1+x_2i_2+ \cdots+x_ni_n+ \bar{x}r$. Then $k_2 = k_1+ sx_1i_1+sx_2i_2+ \cdots+sx_ni_n+ s\bar{x}r$, $a^{k_2} = a^{k_1+ sx_1i_1+sx_2i_2+ \cdots+sx_ni_n+ s\bar{x}r} = a^{k_1+ sx_1i_1+sx_2i_2+ \cdots+sx_ni_n}$. Next, we will construct a weakly connected path joining these two points.
    \par If $x_j > 0$, then the path is 
    \begin{gather*}
        a^{k_1+ sx_1i_1+sx_2i_2+ \cdots+sx_{j-1}i_{j-1}} \to a^{k_1+ sx_1i_1+sx_2i_2+ \cdots+sx_{j-1}i_{j-1} + i_j} \\ \to \cdots \to a^{k_1+ sx_1i_1+sx_2i_2+ \cdots+sx_{j-1}i_{j-1} + sx_ji_j}.
    \end{gather*}
    \par If $x_j < 0$, then the path is 
    \begin{gather*}
        a^{k_1+ sx_1i_1+sx_2i_2+ \cdots+sx_{j-1}i_{j-1}} \gets a^{k_1+ sx_1i_1+sx_2i_2+ \cdots+sx_{j-1}i_{j-1} - i_j} \\ \gets \cdots \gets a^{k_1+ sx_1i_1+sx_2i_2+ \cdots+sx_{j-1}i_{j-1} + sx_ji_j}.
    \end{gather*} 
    For any element $a_i\in \textbf{\textit{NCP}}$, there exists some positive integer multiple $lr$ of its associated period r such that $a^{i+lr}\in \textbf{\textit{CP}}$. Consequently, even if a path involves elements from NCP, we can replace these elements with their counterparts in CP without altering the validity of the path
    \par If $x_j = 0$, then the path will not contain the arc generated by $x_j$.
    Thus the path forms a weakly connected directed path connecting $a^{k_1}$ and $a^{k_2}$ in $\text{Cay}(\left\langle a \right\rangle, C)$.
    \par Suppose that there exists a weakly connected directed path connecting $a^{k_1}$ and $ a^{k_2}$ in $\text{Cay}(\left\langle a \right\rangle, C)$. Let $a^{k_1} \to a^{k_1+l_1} \to \cdots \to a^{k_1+\sum l_t} = a^{k_2}$ be a weakly connected directed path connecting $a^{k_1}$ and $ a^{k_2}$ in $\text{Cay}(\left\langle a \right\rangle, C)$. If there exists a reverse arrow, then the corresponding $l_j$ is negative. Thus $k_1+\sum l_t = k_2 +sr$ for some $s \ge 0$ and $|k_1-k_2| = |\sum l_t-sr|$. Since $d |r$ and $d|l_j$ for any $j$, then $d||\sum l_t-sr|$. Therefore, $d||k_1-k_2|$, $k_1 \equiv k_2 \pmod d$ and so $(a^{k_1}, a^{k_2}) \in \rho_1$.
\end{proof}

\begin{proposition}\label{pro: unique rhoclass}
Suppose that $C = \left\{a^{i_1}, a^{i_2}, \cdots, a^{i_n}\right\} \subseteq \left\langle a \right\rangle$. For any $a^u \in \textbf{\textit{NCP}}$ and $a^{i_p}, a^{i_q} \in C$, if $t_1, t_2$ be the smallest positive integers such that $a^{u+t_1i_p}, a^{u+t_2i_q} \in \textbf{\textit{CP}}$, then $(a^{u+t_1i_p}, a^{u+t_2i_q}) \in \rho_1$.
\end{proposition}
\begin{proof}
    Suppose that $u+t_1i_p-m = s_1r+c_1,0 \le c_1<r$ and $u+t_2i_q-m = s_2r+c_2,0 \le c_2<r$, then $a^{u+t_1i_p} = a^{m+c_1}$ and $a^{u+t_2i_q} = a^{m+c_2}$. Since $d|(t_1i_p-t_2i_q)$ and $d|r$,  $d|(c_1-c_2)$ and so $c_1 \equiv c_2 \pmod d$. Thus $(a^{m+c_1},a^{m+c_2}) \in \rho_1$, moreover, $(a^{u+t_1i_p}, a^{u+t_2i_q}) \in \rho_1$.
\end{proof}

Let $\textbf{\textit{CP}}/\rho_1 = \left\{\bar{C}_1, \bar{C}_2, \cdots, \bar{C}_d\right\}$ and $\textbf{\textit{NCP}}/\rho_2 = \left\{\bar{N}_1, \bar{N}_2, \cdots, \bar{N}_t\right\}(t \le d)$. Without loss of generality, the $\rho_1$ class of $a^m,a^{m+1},\cdots,a^{m+ \gcd(r,d)-1}$ are denoted as $\bar{C}_1, \bar{C}_2, \cdots, \bar{C}_d$ respectively. It follows that $|\bar{C}_1| = |\bar{C}_2| = \cdots= |\bar{C}_d| = \frac{r}{d}$.
\par For arbitrary $\bar{N}_p \in \textbf{\textit{NCP}}/\rho_2$ and $a^u \in \bar{N}_p$, by Proposition~\ref{pro: unique rhoclass}, there is a unique $\rho_1$-class $\bar{C}_q$, such that $a^{u+t_ki_k} \in \bar{C}_q$ for any $a^{i_k} \in C$ and some $t_k \in \mathbb{N}$. Let $\bar{A}_q = \bar{C}_q \cup \bar{N}_p$. Finally, if there exists a $\bar{C}_q$ such that no $\bar{N}_p$ matches it, then let $\bar{A}_q = \bar{C}_q$. By the definition of $\rho_2$, $\left\{\bar{A}_1, \bar{A}_2, \cdots, \bar{A}_d\right\}$ is a partition of $\left\langle a \right\rangle$. By Proposition~\ref{pro: weakly connected}, for arbitrary $\bar{A}_r \in \left\{\bar{A}_1, \bar{A}_2, \cdots, \bar{A}_d\right\}$, the induced subgraph of $\bar{A}_r$ in $\text{Cay}(\left\langle a \right\rangle, C)$ is a weakly connected component of $\text{Cay}(\left\langle a \right\rangle, C)$. We denote such a partition by $Par(C,\rho_1,\rho_2)$ i.e.,$Par(C,\rho_1,\rho_2) = \left\{\bar{A}_1, \bar{A}_2, \cdots, \bar{A}_d\right\}$.

\begin{proposition}\label{pro: number}
Suppose that $C = \left\{a^{i_1}, a^{i_2}, \cdots, a^{i_n}\right\} \subseteq \left\langle a \right\rangle$. Let $d > 1$, for arbitrary $a^u\in \left\langle a \right\rangle \setminus \left\{a^{m-1}\right\}$, if $a^u \in \bar{A}_p$, then $a^{u+1} \in \bar{A}_{p+1}$ (if $p = d, \text{ then } a^{u+1} \in \bar{A}_1$). Further, either 
\begin{align*}
    |\bar{A}_1| = |\bar{A}_2| = \cdots = |\bar{A}_d|
\end{align*}
or there exist $1 \le k_1 \le k_2 \le d$ such that, for all \\$s \in \left\{1, 2, \cdots, k_1-1, k_2+1, \cdots, d\right\}$
\begin{gather*}
    |\bar{A}_{k_1}| = |\bar{A}_{k_1+1}| = \cdots = |\bar{A}_{k_2}| \\
    \text{ and } |\bar{A}_{k_1}|\pm 1 = |\bar{A}_s|.
\end{gather*}
\end{proposition}
\begin{proof}
    Suppose that $\bar{A}_p = \bar{N}_p \cup \bar{C}_p$. If $a^u \in \bar{N}_p$, then there exists $a^{i_k} \in C$, such that $a^{u+t_ki_k}\in \bar{C}_p$ for some $t_k \in \mathbb{N}$. This means $u+t_ki_k \equiv m+p-1 \pmod d$. Then $u+t_ki_k+1 \equiv m+p-1+1 \pmod d$. Since $u \ne m-1$, $a^{u+t_ki_k+1} \in \bar{C}_{p+1}$ and $a^{u+1} \in \bar{N}_{p+1}$. Thus $a^{u+1} \in \bar{A}_{p+1}$. If $p = d$, then $u+t_ki_k+1 \equiv m \pmod d$. Thus $a^{u+t_ki_k+1} \in \bar{C}_{1}$, $a^{u+1} \in \bar{N}_1$. Therefore $a^{u+1} \in \bar{A}_1$.
    \par If $a^u \in \bar{C}_p$, then $u \equiv m+p-1 \pmod d$ and $u+1 \equiv m+p-1+1 \pmod d$. Thus $a^{u+1} \in \bar{C}_{p+1}$, $a^{u+1} \in \bar{A}_{p+1}$. If $p = d$, then $u+1 \equiv m \pmod d$. Thus $a^{u+1} \in \bar{C}_{1}$ and so $a^{u+1} \in \bar{A}_{1}$.
    \par Let $m-1 = td +c \text{ } (0 \le c < d)$. If $c =0$, by $|\bar{C}_1| = |\bar{C}_2| = \cdots= |\bar{C}_d|$ and the above conclusion, $|\bar{A}_1| = |\bar{A}_2| = \cdots = |\bar{A}_d|$. Suppose that $c \ne 0$ and $a \in \bar{A}_{k_1}$. If $k_1+c-1 \le d$, then $a^{td+1} \in \bar{A}_{k_1}, a^{td+2} \in \bar{A}_{k_2}, \cdots, a^{td+c} \in \bar{A}_{k_1+c-1}$. Let $k_2 = k_1+c-1$, then $|\bar{A}_{k_1}| = |\bar{A}_{k_1+1}| = \cdots = |\bar{A}_{k_2}|$ and $|\bar{A}_{k_1}|-1 = |\bar{A}_s|$ for all $s \in \left\{1, 2, \cdots, k_1-1, k_2+1, \cdots, d\right\}$. If $k_1+c-1 > d$, then $a^{td+1} \in \bar{A}_{k_1}, a^{td+2} \in \bar{A}_{k_2}, \cdots, a^{td+c_1} \in \bar{A}_{d}, a^{td+c_1+1} \in \bar{A}_{1}, \cdots, a^{td+c} \in \bar{A}_{k_2}$. By exchanging the values of $k_1$ and $k_2$, we obtain $|\bar{A}_{k_1}|+1 = |\bar{A}_s|$ for all $s \in \left\{1, 2, \cdots, k_1-1, k_2+1, \cdots, d\right\}$.
\end{proof}

\begin{proposition}\label{pro: isomorphism of connected component}
For arbitrary $\bar{A}_p,\bar{A}_q \in \left\{\bar{A}_1, \bar{A}_2, \cdots, \bar{A}_d\right\}$, if $|\bar{A}_p| = |\bar{A}_q|$, then $\text{Cay}(\left\langle a \right\rangle, C)[\bar{A}_p] \simeq \text{Cay}(\left\langle a \right\rangle, C)[\bar{A}_q]$. 
\end{proposition}
\begin{proof}
    Suppose that $|\bar{A}_{k_1}| = |\bar{A}_{k_1+1}|$. By Proposition~\ref{pro: number}, define a mapping $f$ from $\bar{A}_{k_1} $ to $ \bar{A}_{k_1+1} $ as follows: 
    \begin{align*}
        f(a^u) = a^{u+1}.
    \end{align*} 
    It is clear that $f$ is a bijection. Then $(a^u,a^{u+i_k}) \in E(\text{Cay}(\left\langle a \right\rangle, C)[\bar{A}_{k_1}]) \iff (f(a^u), f(a^{u+i_k})) = (a^{u+1}, a^{u+i_k+1}) \in E(\text{Cay}(\left\langle a \right\rangle,\allowbreak C)[\bar{A}_{k_1+1}])$. Thus $|\bar{A}_{k_1}| = |\bar{A}_{k_1+1}| \Longrightarrow \text{Cay}(\left\langle a \right\rangle, C)[\bar{A}_{k_1}] \cong \text{Cay}(\left\langle a \right\rangle,\allowbreak C)[\bar{A}_{k_1+1}]$. This shows that $f$ is a labelled colour-permuting isomorphism. Further, $|\bar{A}_p| = |\bar{A}_q| \Longrightarrow \text{Cay}(\left\langle a \right\rangle, C)[\bar{A}_p] \allowbreak \simeq \text{Cay}(\left\langle a \right\rangle, C)[\bar{A}_q]$.
\end{proof}

\begin{example}
Let $\left\langle a \right\rangle$ be a monogenic semigroup whose index is $4$ and period is $12$, and 
\begin{align*}
    \left\langle a \right\rangle =\left\{a, a^2, \cdots , a^{15}\right\}.
\end{align*}
Let $C = \left\{a^3, a^6, a^9\right\}$. Then
\begin{align*}
    \textbf{\textit{NCP}}=\{a,a^2, a^3 \},
    \textbf{\textit{CP}}=\{a^4,a^5,a^6, a^7, a^8, a^9, a^{10}, a^{11}, a^{12}, a^{13}, a^{14}, a^{15}\}.
\end{align*}
Then
\begin{gather*}
    \textbf{\textit{NCP}}/\rho_2 = \{\{a\}, \{a^2\}, \{a^3\}\},\\
    \textbf{\textit{CP}}/\rho_1 = \{\{a^4, a^7, a^{10}, a^{13}\}, \{a^5, a^8, a^{11}, a^{14}\}, \{a^6, a^9, a^{12}, a^{15}\}\}.
\end{gather*}
For $a^3 \in C$, $a \cdot a^3 = a^4 \in \bar{C}_1$, then $\bar{A}_1 = \{a, a^4, a^7, a^{10}, a^{13}\}$; $a^2 \cdot a^3 = a^5 \in \bar{C}_2$, then $\bar{A}_2 = \{a^2, a^5, a^8, a^{11}, a^{14}\}$; $a^3 \cdot a^3 = a^6 \in \bar{C}_3$, then $\bar{A}_3 = \{a^3, a^6, a^9, a^{12}, a^{15}\}$. In this way, we obtain the following partition of $\left\langle a \right\rangle$
\begin{align*}
    \left\{\left\{a, a^4, a^7, a^{10}, a^{13}\right\}, \left\{a^2, a^5, a^8, a^{11}, a^{14}\right\}, \left\{a^3, a^6, a^9, a^{12}, a^{15}\right\}\right\}.
\end{align*}
The induced subgraph corresponding to each $\bar{A}_i(i= 1,2,3) $ is shown below:
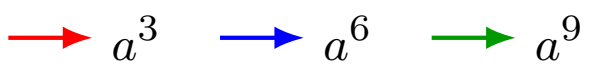
\begin{center}
\makebox[\textwidth][c]{%
\begin{tikzpicture}[scale=0.7,
    vertex/.style={
        circle,
        draw,
        fill=blue!12,
        minimum size=8mm,
        inner sep=0pt
    },
    gen3/.style={-{Latex[length=2mm]}, thick, red},
    gen6/.style={-{Latex[length=2mm]}, thick, blue},
    gen9/.style={-{Latex[length=2mm]}, thick, green!60!black},
    every node/.style={font=\small}
]
\begin{scope}[shift={(-6,0)}]
    \node[vertex] (v1)  at (2.8,2.2)    {$a$};
    \node[vertex] (v4)  at (-0.5,1.2)  {$a^4$};
    \node[vertex] (v7)  at (1.2,1.2)   {$a^7$};
    \node[vertex] (v10) at (1.2,-1.2)  {$a^{10}$};
    \node[vertex] (v13) at (-0.5,-1.2) {$a^{13}$};
    \draw[gen3] (v1)  to[bend right=10] (v4);
    \draw[gen3] (v4)  to[bend left=12] (v7);
    \draw[gen3] (v7)  to[bend left=12] (v10);
    \draw[gen3] (v10) to[bend left=12] (v13);
    \draw[gen3] (v13) to[bend left=12] (v4);
    \draw[gen6] (v1)  -- (v7);
    \draw[gen6] (v4)  to[bend left=12] (v10);
    \draw[gen6] (v10) to[bend left=12] (v4);
    \draw[gen6] (v7)  to[bend left=12] (v13);
    \draw[gen6] (v13) to[bend left=12] (v7);
    \draw[gen9] (v1)  to[bend left=20] (v10);
    \draw[gen9] (v4)  to[bend right=20] (v13);
    \draw[gen9] (v13) to[bend right=20] (v10);
    \draw[gen9] (v10) to[bend right=20] (v7);
    \draw[gen9] (v7)  to[bend right=20] (v4);
\end{scope}
\node at (-5,-2.5) {Figure 5: Cay($\langle a\rangle$, $C$)[$\bar{A}_1$]};
\begin{scope}[shift={(0,0)}]
    \node[vertex] (v2)  at (2.8,2.2)    {$a^2$};
    \node[vertex] (v5)  at (-0.5,1.2)  {$a^5$};
    \node[vertex] (v8)  at (1.2,1.2)   {$a^8$};
    \node[vertex] (v11) at (1.2,-1.2)  {$a^{11}$};
    \node[vertex] (v14) at (-0.5,-1.2) {$a^{14}$};
    \draw[gen3] (v2)  to[bend right=10] (v5);
    \draw[gen3] (v5)  to[bend left=12] (v8);
    \draw[gen3] (v8)  to[bend left=12] (v11);
    \draw[gen3] (v11) to[bend left=12] (v14);
    \draw[gen3] (v14) to[bend left=12] (v5);
    \draw[gen6] (v2)  -- (v8);
    \draw[gen6] (v5)  to[bend left=12] (v11);
    \draw[gen6] (v11) to[bend left=12] (v5);
    \draw[gen6] (v8)  to[bend left=12] (v14);
    \draw[gen6] (v14) to[bend left=12] (v8);
    \draw[gen9] (v2)  to[bend left=20] (v11);
    \draw[gen9] (v5)  to[bend right=20] (v14);
    \draw[gen9] (v14) to[bend right=20] (v11);
    \draw[gen9] (v11) to[bend right=20] (v8);
    \draw[gen9] (v8)  to[bend right=20] (v5);
\end{scope}
\node at (1,-2.5) {Figure 6: Cay($\langle a\rangle$, $C$)[$\bar{A}_2$]};
\begin{scope}[shift={(6,0)}]
    \node[vertex] (v3)  at (2.8,2.2)    {$a^3$};
    \node[vertex] (v6)  at (-0.5,1.2)  {$a^6$};
    \node[vertex] (v9)  at (1.2,1.2)   {$a^9$};
    \node[vertex] (v12) at (1.2,-1.2)  {$a^{12}$};
    \node[vertex] (v15) at (-0.5,-1.2) {$a^{15}$};
    \draw[gen3] (v3)  to[bend right=10] (v6);
    \draw[gen3] (v6)  to[bend left=12] (v9);
    \draw[gen3] (v9)  to[bend left=12] (v12);
    \draw[gen3] (v12) to[bend left=12] (v15);
    \draw[gen3] (v15) to[bend left=12] (v6);
    \draw[gen6] (v3)  -- (v9);
    \draw[gen6] (v6)  to[bend left=12] (v12);
    \draw[gen6] (v12) to[bend left=12] (v6);
    \draw[gen6] (v9)  to[bend left=12] (v15);
    \draw[gen6] (v15) to[bend left=12] (v9);
    \draw[gen9] (v3)  to[bend left=20] (v12);
    \draw[gen9] (v6)  to[bend right=20] (v15);
    \draw[gen9] (v15) to[bend right=20] (v12);
    \draw[gen9] (v12) to[bend right=20] (v9);
    \draw[gen9] (v9)  to[bend right=20] (v6);
\end{scope}
\node at (7,-2.5) {Figure 7: Cay($\langle a\rangle$, $C$)[$\bar{A}_3$]};
\begin{scope}[shift={(-1.7,-3.8)}]
    \draw[gen3] (0,0) -- +(0.8,0) node[right,black] {$a^3$};
    \draw[gen6] (2,0) -- +(0.8,0) node[right,black] {$a^6$};
    \draw[gen9] (4,0) -- +(0.8,0) node[right,black] {$a^9$};
\end{scope}

\end{tikzpicture}%
}
\end{center}
\end{example}

\section{The labelled colour-permuting isomorphisms of Cayley graphs}

Having determined the weakly connected components in Section 3, we now turn to the isomorphism problem. The key idea is to compare connection sets through the decomposition into non-cyclic and cyclic parts, and to encode the componentwise behaviour by admissible ordered-reduced connection sets. This leads to a complete classification of labelled colour-permuting isomorphisms in the cases considered below.
\par For the convenience of the following discussion, the connection set $C$ is divided into the following three cases: (i) $C \subseteq \textbf{\textit{NCP}}$, (ii) $C\subseteq \textbf{\textit{CP}}$ and (iii) both $C \cap \textbf{\textit{NCP}} \ne \emptyset \text{ and } C\cap \textbf{\textit{CP}} \ne \emptyset$.
\par It is natural to ask how far the CI problem can be extended from groups to semigroups. Compared with the group case, the absence of inverses and the possible presence of non-periodic elements make the semigroup case substantially different. In particular, graph isomorphisms need not arise from algebraic automorphisms in the same way as in the classical group setting.

\begin{proposition}\label{pro: NCP}
Let $\left\langle a \right\rangle$ be a monogenic semigroup.
Suppose that $C_1 = \left\{a^{i_1}, a^{i_2}, \cdots, a^{i_n}\right\}$ and $i_1,i_2,...,i_k< m-1,i_{k+1}, i_{k+2}, \cdots ,i_n \ge m-1$. If there exists $C_2 \subseteq \left\langle a \right\rangle$, such that $\text{Cay}(\left\langle a \right\rangle,C_1) \simeq \text{Cay}(\left\langle a \right\rangle,C_2)$, then $a^{i_t} \in C_2$ for all $t \in \left\{1,2, \cdots,k\right\}$.
\end{proposition}
\begin{proof}
    In $\text{Cay}(\left\langle a \right\rangle, C_1)$, the in-degree of vertices $a, a^2,\cdots,a^{i_1}$ are all, the in-degree of vertices $a^{i_1+1}, a^{i_1+2},\cdots, a^{i_2}$ are all $1$, and so on. The in-degree of vertices $a^{i_{k-1}+1}, \allowbreak a^{i_{k-1}+2}, \cdots,a^{i_k}$ are all $k-1$. And the in-degree of  vertices $a^{i_k+1}, a^{i_k+2},\cdots,a^{m+r-1}$ are more then $k-1$.
    \par Since $\text{Cay}(\left\langle a \right\rangle,C_1) \simeq \text{Cay}(\left\langle a \right\rangle,C_2)$, the number of vertices with $0$ in-degree in $\text{Cay}(\left\langle a \right\rangle,C_2)$ is $i_1$. Then $a^{i_1} \in C_2$. Similarly, $a^{i_2}, \cdots,a^{i_k} \in C_2$. 
\end{proof}

\begin{corollary}\label{coro: 2}
Let $C_1 = \left\{a^{i_1}, a^{i_2}, \cdots, a^{i_n}\right\} \subseteq \textbf{NCP}$. If there exists $C_2 \subseteq \left\langle a \right\rangle$, such that $\text{Cay}(\left\langle a \right\rangle,C_1) \simeq \text{Cay}(\left\langle a \right\rangle,C_2)$, then $C_1 = C_2$. Further, suppose that $a^i, a^j \in \textbf{\textit{NCP}}$, $\text{Cay}(\left\langle a \right\rangle,\{a^i\}) \simeq  \text{Cay}(\left\langle a \right\rangle,\{a^j\})$ if and only if $i=j$.
\end{corollary}
\begin{proof}
    Since $C_1 = \left\{a^{i_1}, a^{i_2}, \cdots, a^{i_n}\right\} \subseteq \textbf{NCP}$, by Proposition \ref{pro: NCP}, it is clear that $C_1 = C_2$.
\end{proof}

We put an order on the reduced connection sets thus obtained. First list the elements of $\textbf{\textit{NCP}}\setminus\{a^{m-1}\}$ in increasing order of their exponents. Then list the elements of $\textbf{\textit{CP}}$ in increasing order of $\gcd(i,r)$.
\par When ordering the connection sets $C_1$ and $C_2$, if two elements have the same greatest common divisor of their exponents with $r$, then the order of these two elements is not required. This implies that the reduced connection sets resulting from the ordering are not uniquely determined. For example, when $r=8$, both
\begin{align*}
    (a^3,a^2,a^6)
    \quad\text{and}\quad
    (a^3,a^6,a^2)
\end{align*}
are reduced connection sets that satisfy the ordering because of $\gcd(2,8)=\gcd(6,8)=2$. 
\par We call such an ordered reduced connection set an admissible ordered-reduced connection set.

\begin{proposition}\label{pro: weakly connected component to another Cay}
    Let $\left\langle a \right\rangle$ be a monogenic semigroup, $C = \left\{a^{i_1}, a^{i_2}, \cdots, a^{i_n}\right\}\allowbreak \subseteq \left\langle a \right\rangle$ be an admissible ordered-reduced connection set and $\gcd(C, r) = \gcd(i_1, r)$. If $\text{Cay}(\left\langle a \right\rangle,C)$ has more than one weakly connected component, then each weakly connected component is labelled colour-permuting isomorphic to the Cayley graph of a monogenic semigroup.
\end{proposition}
\begin{proof}
    Suppose that $\bar{A}_k$ is a weakly connected component of $\text{Cay}(\left\langle a \right\rangle,C)$. Since $\gcd(C, r) = \gcd(i_1, r) = d$, the elements of $\bar{A}_k$ can be denoted as follows,
    \begin{align*}
        \bar{A}_k = \{a^u, a^{u+d}, \cdots, a^{u+(m'-2)d}, a^{u+(m'-1)d}, \cdots, a^{u+(m'+r'-2)d}\},
    \end{align*}
    where $m' = |\bar{A}_k| - |\bar{C}_k| +1$ and $r' = \frac{r}{d}$.
    Let
    \begin{align*}
        \left\langle b \right\rangle = \{b, b^2, \cdots, b^{m'-1}, b^{m'}, \cdots, b^{m'+r'-1}\}
    \end{align*}
    be a monogenic semigroup with index $m'$ and period $r'$,
    \begin{align*}
        C' = \{b^{\frac{i_1}{d}}, b^{\frac{i_2}{d}}, \cdots, b^{\frac{i_n}{d}}\}.
    \end{align*}
    Since $m'+r'-1 \ge \frac{i_n}{d} $, $C' \subseteq \left\langle b \right\rangle$. Consider the Cayley graph $\text{Cay}(\left\langle b \right\rangle,C')$.
    Let $f$ be a mapping from $\bar{A}_k$ to $\left\langle b \right\rangle$ and $f(a^{u+td}) = b^{t+1}$ for any $a^{u+td} \in \bar{A}_k$. It is clear that $f$ is a bijection. Let $\varphi$ be a bijection from $C$ to $C'$ and $\varphi(a^{i_k}) = b^{\frac{i_k}{d}}$. For any $a^{i_k} \in C$ and $a^{u+td} \in \left\langle a \right\rangle$, 
    \begin{align*}
        f(a^{i_k} \cdot a^{u+td}) = f(a^{u+(t+\frac{i_k}{d})d}) = b^{t+1+\frac{i_k}{d}} = b^{\frac{i_k}{d}}\cdot b^{t+1} = \varphi(a^{i_k})\cdot f(a^{u+td}).
    \end{align*}
    Thus $\bar{A}_k$ is labelled colour-permuting isomorphic to $\text{Cay}(\left\langle b \right\rangle,C')$.
\end{proof}

\begin{proposition}\label{pro: gcd}
Let $C_1, C_2 \subseteq \textbf{\textit{CP}}$, $\gcd(C_1, r) = \gcd(C_2,r) = \gcd(i_\alpha, r) = \gcd(j_\beta,r)=1$ for some $a^{i_\alpha} \in C_1$ and $a^{j_\beta} \in C_2$. If $\text{Cay}(\left\langle a \right\rangle,C_1) \simeq \text{Cay}(\left\langle a \right\rangle,C_2)$, then there are admissible ordered-reduced connection sets of $C_1$ and $C_2$, and with respect to this ordering they can be expressed as:
\begin{align*}
    C_1 = \left\{a^{i_1}, a^{i_2}, \cdots, a^{i_n}\right\} \text{ and } C_2=\left\{a^{j_1}, a^{j_2}, \cdots, a^{j_n}\right\},
\end{align*}
such that for all $k \in \left\{1,2,\cdots, n\right\}$
    \begin{gather*}
        \gcd(i_k,r) = \gcd(j_k,r),
         \frac{i_k}{i_1} \equiv \frac{j_k}{j_1} \pmod r.
    \end{gather*}
\end{proposition}
\begin{proof}
    Suppose that $C_1 = \left\{a^{i_1}, a^{i_2}, \cdots, a^{i_n}\right\}$ be an admissible ordered-reduced connection set and with respect to this ordering they can be expressed as $C_1 = \left\{a^{i_1}, a^{i_2}, \cdots, a^{i_n}\right\}$. Let $f$ be the labelled colour-permuting isomorphism, $\varphi$ be the bijection from $C_1$ to $C_2$, $C_2 = \left\{a^{j_1}, a^{j_2}, \cdots, a^{j_n}\right\}$ and $\varphi(a^{i_k}) = a^{j_k}$. For arbitrary $a^{i_k} \in C_1$, 
    \begin{align*}
        a^m \to a^{m+i_k} \to \cdots \to a^{m+(\frac{r}{\gcd(i_k),r}-1)i_k} \to a^{m+\frac{r}{\gcd(i_k,r)}i_k} = a^m
    \end{align*}
    is a directed circle in $\text{Cay}(\left\langle a \right\rangle,C_1)$. Hence 
    \begin{align*}
        f(a^m) \to f(a^{m+i_k}) \to \cdots \to f(a^{m+(\frac{r}{\gcd(i_k),r}-1)i_k}) \to f(a^{m+\frac{r}{\gcd(i_k,r)}i_k}) = f(a^m)
    \end{align*} 
    is a directed circle in $\text{Cay}(\left\langle a \right\rangle,C_2)$. Since $\text{Cay}(\left\langle a \right\rangle,C_1) \simeq \text{Cay}(\left\langle a \right\rangle,C_2)$, 
    \begin{gather*}
        f(a^{m+i_k}) = \varphi(a^{i_k})f(a^m), f(a^{m+2i_k}) = \varphi(a^{i_k})\varphi(a^{i_k})f(a^m), \cdots,\\
        f(a^{m+(\frac{r}{\gcd(i_k),r}-1)i_k}) = (\varphi(a^{i_k}))^{\frac{r}{\gcd(i_k, r)}-1}f(a^m).
    \end{gather*}
    Thus $\varphi(a^{i_k}) \in C_2$ and this directed circle 
    \begin{align*}
        f(a^m) \to f(a^{m+i_k}) \to \cdots \to f(a^{m+(\frac{r}{\gcd(i_k),r}-1)i_k}) \to f(a^{m+\frac{r}{\gcd(i_k,r)}i_k}) = f(a^m)
    \end{align*} starts at $f(a^m)$ and is generated by $\varphi(a^{i_k}) = a^{j_k}$. Then the length of this directed cycle is $\frac{r}{\gcd(j_k, r)}$. This means that $\gcd(i_k,r) = \gcd(j_k,r)$ for all $k \in \left\{1,2,\cdots, n\right\}$. And it is clear that $C_2 = \left\{a^{j_1}, a^{j_2}, \cdots, a^{j_n}\right\}$ is an admissible ordered-reduced connection set of $C_2$.
    \par Take $a^m \in \textbf{\textit{CP}}$ and suppose that $f(a^m) = a^{m+sj_1}$ for some $s \in \mathbb{N}$. Since both $C_1 = \left\{a^{i_1}, a^{i_2}, \cdots, a^{i_n}\right\} \text{ and } C_2=\left\{a^{j_1}, a^{j_2}, \cdots, a^{j_n}\right\}$ are the admissible ordered-reduced connection sets, it follows that that $\gcd(C_1, r) = \gcd(C_2,r) = \gcd(i_1, r) = \gcd(j_1,r)=1$ and directed circles 
    \begin{align*}
        a^{m} \to a^{m+i_1} \to \cdots \to a^{m+(\frac{r}{\gcd(i_1,r)}-1)i_1} \to a^{m}
    \end{align*}
    and 
    \begin{align*}
        a^{m+t_2j_1} \to a^{m+t_2j_1+j_1} \to \cdots \to a^{m+t_2j_1+(\frac{r}{\gcd(j_1,r)}-1)j_1} \to a^{m+t_2j_1}
    \end{align*}
    are the longest circles in $\text{Cay}(\left\langle a \right\rangle,C_1) \text{ and } \text{Cay}(\left\langle a \right\rangle,C_2)$. 
    \par For all $k \in \left\{1,2,\cdots, n\right\}$, there exist $t_{k_1}, t_{k_2} \in \mathbb{N}$, such that 
    \begin{gather*}
        a^{m+t_{k_1}i_1} = a^{i_k} \cdot a^m,\\
        a^{m+sj_1+t_{k_2}j_1} =a^{j_k} \cdot a^{m+sj_1}.
    \end{gather*}
    And we have that 
    \begin{gather*}
        f(a^{i_k}\cdot a^m) = \varphi(a^{i_k})f(a^m) = a^{j_k}\cdot a^{m+sj_1},\\
        f(a^{i_k}\cdot a^m) = f(a^{m+t_{k_1}i_1}) = a^{m+sj_1+t_{k_1}j_1} .
    \end{gather*}
    Thus 
    \begin{gather*}
        a^{m+sj_1+t_{k_2}j_1} = a^{m+sj_1+t_{k_1}j_1} 
    \end{gather*}
    and $t_{k_1} \equiv t_{k_2} \pmod r$.
    Therefore
    \begin{gather*}
        t_{k_1} \equiv \frac{i_k}{i_1} \pmod r,
        t_{k_2} \equiv \frac{j_k}{j_1} \pmod r.
    \end{gather*}
    Thus 
    \begin{align*}
        \frac{i_k}{i_1} \equiv \frac{j_k}{j_1} \pmod r
    \end{align*}
    for all $k \in \left\{1,2,\cdots, n\right\}$.
\end{proof}

\begin{lemma}\label{lem: number}
\textup{(\cite{Shoup2009})}
    If $a,b \in \mathbb{N}$ and $a,b > 0$, then $\gcd(a,b) = \gcd(a, (b\bmod a))$.
\end{lemma}

\begin{corollary}\label{coro: same connected components}
Suppose that $C_1, C_2 \subseteq \textbf{\textit{CP}}$ and there exist admissible ordered-reduced connection sets of $C_1$ and $C_2$, and with respect to this ordering they can be expressed as:
\begin{align*}
    C_1 = \left\{a^{i_1}, a^{i_2}, \cdots, a^{i_n}\right\}, C_2=\left\{a^{j_1}, a^{j_2}, \cdots, a^{j_n}\right\}.
\end{align*}
Let $\gcd(C_1, r) = \gcd(C_2,r) = \gcd(i_1, r) = \gcd(j_1,r)=1$. If $\frac{i_k}{i_1} \equiv \frac{j_k}{j_1} \pmod r$ for all $k \in \left\{1,2,\cdots, n\right\}$, then $\gcd(i_k,r) = \gcd(j_k,r)$ for all $k \in \left\{1,2,\cdots, n\right\}$. Furthermore, $\gcd(C_1,r) = \gcd(C_2, r)$ and $Par(C_1, \rho_1,\rho_2) = Par(C_2, \rho_1, \rho_2)$.
\end{corollary}
\begin{proof}
    From the proof of \ref{pro: gcd}, it can be concluded that $t_{k_1}i_1 \equiv i_k \pmod{r}$ and $t_{k_2}j_1 \equiv j_k \pmod{r}$. By Lemma \ref{lem: number}, 
    \begin{gather*}
        \gcd(i_k, r) = \gcd(t_{k_1}i_1,r),
        \gcd(j_k, r) = \gcd(t_{k_2}j_1,r).
    \end{gather*}
    Since $\gcd(i_1, r) = \gcd(j_1,r) = 1$, 
    \begin{gather*}
        \gcd(t_{k_1}i_1,r) = \gcd(t_{k_1}, r),
        \gcd(t_{k_2}j_1,r) = \gcd(t_{k_2}, r).
    \end{gather*}
    Since $t_{k_1} \equiv t_{k_2} \pmod r$, we may assume without loss of generality that $t_{k_1} - t_{k_2} = qr$. So $t_{k_1} = qr + t_{k_2}$. Suppose that $\gcd(t_{k_1}, r) = d$, then $d|t_{k_1}$. Since $d|r$, $d|t_{k_2}$. Thus $d \le \gcd(t_{k_2}, r)$. Suppose that $d' = \gcd(t_{k_2},r)$. Similarly, we can obtain $d' \le \gcd(t_{k_1},r)$. That means $\gcd(t_{k_1},r) \le \gcd(t_{k_2},r)$ and $\gcd(t_{k_2},r) \le \gcd(t_{k_1},r)$. Thus $\gcd(t_{k_1},r) = \gcd(t_{k_2},r)$ and $\gcd(i_k, r) = \gcd(j_k, r)$ for all $k \in \left\{1,2,\cdots, n\right\}$. It follows that $\gcd(C_1,r) = \gcd(C_2, r)$ and $Par(C_1, \rho_1,\rho_2) = Par(C_2, \rho_1, \rho_2)$.
\end{proof}

\begin{lemma}\label{lem: rotating}
Let $a \to a^2 \to \cdots \to a^t \to a$ be a directed circle and suppose that a directed circle $a^{i_1} \to a^{i_2} \to \cdots \to a^{i_t} \to a^{i_1}$ is isomorphic to it. Then this directed circle is obtained from $a \to a^2 \to \cdots \to a^t \to a$ by a rotating, or  by a reflection followed by a rotation.
\end{lemma}
\begin{proof}
    Let $f$ be the graph isomorphism and $f(a) = a^{i_k}$. Then $f(a^2) = a^{i_{k+1}}$ and $f(a^t) = a^{i_{k-1}}$. By analogy, it can be inferred that $a^{i_1} \to a^{i_2} \to \cdots \to a^{i_t} \to a^{i_1}$ is the result of $a \to a^2 \to \cdots \to a^t \to a$ rotating around the center. If the other directed circle is $a^{i_1} \gets a^{i_2} \gets \cdots \gets a^{i_t} \gets a^{i_1}$, it is obtained by flipping $a \to a^2 \to \cdots \to a^t \to a$ first and then rotating it around the center. 
\end{proof}

\begin{remark}\label{rem: 3}
    In the proof of the following theorem, we shall specify the orientation of each directed circle, so that no flipping will be involved.
\end{remark}

\begin{definition}\label{def: EF}
Suppose that $C_1, C_2 \subseteq \textbf{\textit{CP}}$ and there exist admissible ordered-reduced connection sets of $C_1$ and $C_2$, and with respect to this ordering they can be expressed as:
\begin{align*}
    C_1 = \left\{a^{i_1}, a^{i_2}, \cdots, a^{i_n}\right\}, C_2=\left\{a^{j_1}, a^{j_2}, \cdots, a^{j_n}\right\}.
\end{align*}
Let $m-1 = qr +c \text{ }(0 \le c <r)$ and take the last $c$ elements in $\textbf{\textit{NCP}}$. Define the vector $E_k \in \mathbb{N}^c$ denote the arcs that are generated by $a^{i_k}$ and start from these elements,
\begin{gather*}
    E_k = (m+(qr+1+i_k-m)\bmod r, m+(qr+2+i_k-m)\bmod r,  \cdots,\\ 
    m+(m-1+i_k-m)\bmod r ). \\
\end{gather*} 
Similarly, Define the vector $F_k \in \mathbb{N}^c$ denote the arcs that are generated by $a^{j_k}$ and start from these elements,
\begin{gather*}
    F_k = (m+(qr+1+j_k-m)\bmod r, m+(qr+2+j_k-m)\bmod r,  \cdots,\\ 
    m+(m-1+j_k-m)\bmod r ).
\end{gather*} 
\end{definition}

\begin{definition}\label{def: XY}
Suppose that $C_1, C_2 \subseteq \textbf{\textit{CP}}$ and there exist admissible ordered-reduced connection sets of $C_1$ and $C_2$, and with respect to this ordering they can be expressed as:
\begin{align*}
    C_1 = \left\{a^{i_1}, a^{i_2}, \cdots, a^{i_n}\right\}, C_2=\left\{a^{j_1}, a^{j_2}, \cdots, a^{j_n}\right\}.
\end{align*}
Let $\gcd(C_1, r) = \gcd(C_2, r)= \gcd(i_1, r) = \gcd(j_1, r)$. Define the vector $X_k$ to represent the order in which the elements of the set
\begin{align*}
    \{a^{qr+1+i_k}, a^{qr+2+i_k}, \cdots, a^{m-1+i_k}\}
\end{align*}
appear successively on the directed cycle
    \begin{align*}
        a^m \to a^{m+i_1} \to \cdots \to a^{m+(\frac{r}{\gcd(i_1,r)}-1)i_1} \to a^{m}
    \end{align*}
and the vector $Y_k$ to represent the order in which the elements of the set
\begin{align*}
    \left\{a^{qr+1+j_k}, a^{qr+2+j_k}, \cdots, a^{m-1+j_k}\right\}
\end{align*}
appear successively on the directed cycle
    \begin{align*}
        a^m \to a^{m+j_1} \to \cdots \to a^{m+(\frac{r}{\gcd(j_1,r)}-1)j_1} \to a^{m}.
    \end{align*}
\begin{flalign*}
    \quad \quad &X_k, Y_k \in \mathbb{N}^{\lambda}, &\\
    &\lambda=\frac{r}{\gcd(C_1, r)} = \frac{r}{\gcd(C_2, r)}. &\\
    &\text{and let } \pi_i \text{ be canonical morphism.}&\\
    \quad \quad &\pi_t(X_k) = 
    \begin{cases}
        0, & \text{if } m+(ti_1\bmod r) \ne \pi_p(E_k)\\
        &\text{ for any }1 \le p \le c,\\
        1, & \text{if } m+(ti_1\bmod r) = \pi_p(E_k)\\
        &\text{ for some }1 \le p \le c. 
    \end{cases} &\\
    \quad \quad &\pi_t(Y_k) = 
    \begin{cases}
        0, & \text{if } m+(tj_1\bmod r) \ne \pi_q(F_k)\\
        &\text{ for any }1 \le q \le c,\\
        1, & \text{if } m+(tj_1\bmod r) = \pi_q(F_k)\\
        &\text{ for some }1 \le q \le c. 
    \end{cases} &\\
\end{flalign*}

\end{definition}

\begin{theorem}\label{theorem: main}
Let $C_1, C_2 \subseteq \textbf{\textit{CP}}$, $\gcd(C_1, r) = \gcd(C_2,r) = \gcd(i_\alpha, r) = \gcd(j_\beta,r)=1$ for some $a^{i_\alpha} \in C_1$ and $a^{j_\beta} \in C_2$. Then $\text{Cay}(\left\langle a \right\rangle,C_1) \simeq \text{Cay}(\left\langle a \right\rangle,C_2)$ if and only if there exist admissible ordered-reduced connection sets of $C_1$ and $C_2$, and with respect to this ordering they can be expressed as:
\begin{align*}
    C_1 = \left\{a^{i_1}, a^{i_2}, \cdots, a^{i_n}\right\}, C_2=\left\{a^{j_1}, a^{j_2}, \cdots, a^{j_n}\right\}
\end{align*}
satisfy the following conditions.
\begin{flalign*}
    &(1) \text{ For all } k \in \left\{1,2,\cdots, n\right\}, \frac{i_k}{i_1} \equiv \frac{j_k}{j_1} \pmod r.&\\
    &(2) \text{ There exists a rotation } \phi \text{ such that for all } k \in \{1,2,\cdots,n\}, \phi(X_k) = Y_k.
\end{flalign*}
\end{theorem}
\begin{proof}
    \textbf{Sufficiency.} Assume that  there exist admissible ordered-reduced connection sets of $C_1$ and $C_2$, and with respect to this ordering they can be expressed as:
\begin{align*}
    C_1 = \left\{a^{i_1}, a^{i_2}, \cdots, a^{i_n}\right\}, C_2=\left\{a^{j_1}, a^{j_2}, \cdots, a^{j_n}\right\}
\end{align*}
and the two conditions are satisfied. Since $\gcd(i_1, i_2, \cdots,\allowbreak i_n, r) = \gcd(i_1, r) =1$, this graph has only one weakly connected component. Define $\varphi: C_1 \to C_2$, $\varphi(a^{i_k}) = a^{j_k}$ for all $k \in \{1,2, \cdots, n\}$.
    \par Suppose that $\phi$ maps the first component of the preimage to the $x^{th}$ component of the image. Define $f_0: \textbf{\textit{CP}} \to \textbf{\textit{CP}}$, $f_0(a^{m+ti_1}) = a^{m+(x+t)j_1}$ for all $a^{m+ti_1} \in \textbf{\textit{CP}}$. It is clear that $f_0$ is a bijection. For all $k \in \left\{1,2,\cdots, n\right\}$, there exist $t_{k_1}, t_{k_2} \in \mathbb{N}$, such that 
    \begin{gather*}
        a^{m+t_{k_1}i_1} = a^{i_k} \cdot a^m,\\
        a^{m+xj_1+t_{k_2}j_1} =a^{j_k} \cdot a^{m+xj_1}.
    \end{gather*}
    It is clear that 
    \begin{gather*}
        i_1t_{k_1} \equiv i_k \pmod r,
        j_1t_{k_2} \equiv j_k \pmod r.
    \end{gather*}
    Since 
    \begin{align*}
        \frac{i_k}{i_1} \equiv \frac{j_k}{j_1} \pmod r
    \end{align*}
    for all $k \in \left\{1,2,\cdots, n\right\}$. Then $i_k \equiv \frac{i_1j_k}{j_1} \pmod r$ and $\frac{j_k}{j_1} \equiv t_{k_1} \pmod r$. Thus $j_1t_{k_1} \equiv j_1t_{k_2} \pmod r$. 
    \par For any $a^{m+ti_1} \in \textbf{\textit{CP}}$ and $a^{i_k} \in C_1$, 
    \begin{gather*}
        f_0(a^{i_k} \cdot a^{m+ti_1}) = f_0(a^{t_{k_1}i_1} \cdot a^{m+ti_1}) = f_0(a^{m+(t+t_{k_1})i_1}) = a^{m+(t+t_{k_1}+x)j_1} = \\
        a^{m+(t+x)j_1+t_{k_2}j_1} = a^{j_k}\cdot a^{m+(x+t)j_1} = \varphi(a^{i_k}) \cdot f_0(a^{m+ti_1})
    \end{gather*}
    \par Let $m-1 = qr +c \text{ }(0 \le c <r)$. For any $a^u \in \{a, a^2, \cdots, a^r\}$, $a^{i_k}\cdot a^u \in \textbf{\textit{CP}}$ for any $a^{i_k} \in C_1$. Since $\phi(X_k) = Y_k$ and for all $k \in \left\{1,2,\cdots, n\right\}$
    \begin{align*}
        \frac{i_k}{i_1} \equiv \frac{j_k}{j_1} \pmod r.
    \end{align*}
    there exists a unique element $a^v \in \{a, a^2, \cdots, a^r\}$, such that $a^{j_k}\cdot a^v = f_0(a^{i_k}\cdot a^u)$. Let $f_1(a^u) = a^v$, it follows that $f_1$ is a bijection from $\{a, a^2, \cdots, a^r\}$ to $\{a, a^2, \cdots, a^r\}$. Similarly we obtain bijections $f_2, f_3, \cdots, f_q$. If $q = 0$, then we do not need to define these maps.
    \par For any $a^u \in \{a^{qr+1}, a^{qr+2}, \cdots, a^{m-1}\}$, $a^{i_k}\cdot a^u \in \textbf{\textit{CP}}$ for any $a^{i_k} \in C_1$. Since $\phi(X_k) = Y_k$, there exists a unique $a^v \in \{a^{qr+1}, a^{qr+2}, \cdots, a^{m-1}\}$, such that $a^{j_k}\cdot a^v = f_0(a^{i_k}\cdot a^u)$. Let $f_{q+1}(a^u) = a^v$, we have that $f_{q+1}$ is a bijection.
    \par Define $f: \left\langle a \right\rangle \to \left\langle a \right\rangle$, 
    \begin{gather*}
    f(a^u)=
        \begin{cases}
        f_0(a^u), \text{ if } a^u \text{ in } \textbf{\textit{CP}},\\
        f_1(a^u), \text{ if } a^u \in \{a, a^2, \cdots, a^r\},\\
        \quad \quad \quad\quad\quad\quad \cdots\cdots\\
        f_q(a^u), \text{ if } a^u \in \{a^{(q-1)r+1}, a^{(q-1)r+2}, \cdots, a^{qr}\},\\
        f_{q+1}(a^u), \text{ if } a^u \in \{a^{qr+1}, a^{qr+2}, \cdots, a^{m-1}\}.
    \end{cases}
    \end{gather*}
    Then for any $a^{u} \in \left\langle a \right\rangle$ and $a^{i_k} \in C_1$, $f(a^{i_k}\cdot a^u) = \varphi(a^{i_k})\cdot f(a^u)$. Thus $\text{Cay}(\left\langle a \right\rangle,C_1) \simeq \text{Cay}(\left\langle a \right\rangle,C_2)$.
    \par \textbf{Necessity.} Suppose that there exist admissible ordered-reduced connection sets of $C_1$ and $C_2$, and with respect to this ordering they can be expressed as:
    \begin{align*}
        C_1 = \left\{a^{i_1}, a^{i_2}, \cdots, a^{i_n}\right\}, C_2=\left\{a^{j_1}, a^{j_2}, \cdots, a^{j_n}\right\}
    \end{align*}
    such that for all $k \in \left\{1,2,\cdots, n\right\}$
    \begin{align*}
        \frac{i_k}{i_1} \equiv \frac{j_k}{j_1} \pmod r.
    \end{align*}
    Let $f$ be the labelled colour-permuting isomorphism, $\varphi$ be the bijection from $C_1$ to $C_2$. Since $\gcd(C_1, r) = \gcd(C_2,r) = \gcd(i_1, r) = \gcd(j_1,r)=1$, the circles
    \begin{align*}
        a^{m} \to a^{m+i_1} \to \cdots \to a^{m+(\frac{r}{\gcd(i_1,r)}-1)i_1} \to a^{m},\\
        a^{m} \to a^{m+j_1} \to \cdots \to a^{m+(\frac{r}{\gcd(j_1,r)}-1)j_1} \to a^{m}
    \end{align*}
    are the longest circles in $\text{Cay}(\left\langle a \right\rangle,C_1) \text{ and } \text{Cay}(\left\langle a \right\rangle,C_2)$. By Lemma \ref{lem: rotating} and the definition of $E_k$ and $F_k$, there exists a rotation $\phi$ such that $\phi(E_k) = F_k$ for all $k \in \{1, 2, ,\cdots, n\}$.
\end{proof}

By the definition of $ X_k$ and $Y_k$, we can derive some congruence equations, through which we can readily obtain $ X_k $ and $ Y_k$. Next, we provide the congruence equations for $X_k$, and the equations for $Y_k$ are obtained similarly. $\pi_x(E_k) = m+(qr + x +i_k -m) \mod r$, then for all $k\in \left\{1,2,\cdots , n\right\}$, 
\begin{gather*}
        (1-1)i_1-qr-i_k+m \equiv x \pmod r \\
        (2-1)i_1-qr-i_k+m \equiv x \pmod r \\
        \cdots \\
        (\frac{r}{\gcd(i_1, r)}-1)i_1-qr-i_k+m \equiv x \pmod r
\end{gather*}
If the $t^{th}$ equation has no solution for $1 \le x \le c$, then $\pi_t(E_k)=0$, if the $t^{th}$ equation has solution for $1 \le x \le c$, then $\pi_t(E_k)=1$.

\vspace{1 em}
In \Cref{theorem: main}, we imposed $\gcd(C_1, r) = \gcd(C_2,r) = \gcd(i_\alpha, r) = \gcd(j_\beta,r) = 1$ for some $a^{i_\alpha} \in C_1$ and $a^{j_\beta} \in C_2$. Next, we assume $\gcd(C_1, r) = \gcd(C_2,r) = \gcd(i_\alpha, r) = \gcd(j_\beta,r)$ for some $a^{i_\alpha} \in C_1$ and $a^{j_\beta} \in C_2$.

\begin{theorem}\label{the: many}
    Let $C_1, C_2 \subseteq \textbf{\textit{CP}}$, $\gcd(C_1, r) = \gcd(C_2,r) = \gcd(i_\alpha, r) = \gcd(j_\beta,r)$ for some $a^{i_\alpha} \in C_1$ and $a^{j_\beta} \in C_2$. Suppose that $a^{i_1} \in \bar{A}_{k_1} \in Par(C_1, \rho_1, \rho_2), a^{j_1} \in \bar{A}_{l_1} \in Par(C_2, \rho_1, \rho_2)$. Then $\text{Cay}(\left\langle a \right\rangle,C_1) \simeq \text{Cay}(\left\langle a \right\rangle,C_2)$ if and only if there exist admissible ordered-reduced connection sets of $C_1$ and $C_2$, and with respect to this ordering they can be expressed as:
\begin{align*}
    C_1 = \left\{a^{i_1}, a^{i_2}, \cdots, a^{i_n}\right\}, C_2=\left\{a^{j_1}, a^{j_2}, \cdots, a^{j_n}\right\}
\end{align*}
such that 
\begin{gather*}
    \text{Cay}(\left\langle a \right\rangle, C_1)[\bar{A}_{k_1}] \simeq \text{Cay}(\left\langle a \right\rangle, C_2)[\bar{A}_{l_1}],\\ \text{Cay}(\left\langle a \right\rangle,\allowbreak C_1)[\bar{A}_{k_1+1}] \simeq \text{Cay}(\left\langle a \right\rangle, C_2)[\bar{A}_{l_1+1}].
\end{gather*}
\end{theorem}
\begin{proof}
    Without loss of generality, assume that this Cayley graph has two distinct weakly connected components up to isomorphism. 
    \par Suppose that $\text{Cay}(\left\langle a \right\rangle,C_1) \simeq \text{Cay}(\left\langle a \right\rangle,C_2)$, it is clear that $f$ be a labelled colour-permuting isomorphism and $f$ maps one weakly connected component to another weakly connected component. By Proposition \ref{pro: number}, we can obtain $|\bar{A}_{k_1}| = |\bar{A}_{l_1}|$ and $|\bar{A}_{k_1+1}| = |\bar{A}_{l_1+1}|$. Assume that $f(\bar{A}_{k_1}) = \bar{A}_u$, then $\text{Cay}(\left\langle a \right\rangle, C_1)[\bar{A}_{k_1}] \simeq \text{Cay}(\left\langle a \right\rangle, C_2)[\bar{A}_u]$ and $|\bar{A}_{k_1}| = |\bar{A}_u| = |\bar{A}_{l_1}|$. By Proposition \ref{pro: isomorphism of connected component}, $\text{Cay}(\left\langle a \right\rangle, C_2)\allowbreak [\bar{A}_{l_1}] \simeq \text{Cay}(\left\langle a \right\rangle, C_2)[\bar{A}_u]$. Thus $\text{Cay}(\left\langle a \right\rangle, C_1)[\bar{A}_{k_1}] \simeq \text{Cay}(\left\langle a \right\rangle, C_2)[\bar{A}_{l_1}]$. Similarly, $\text{Cay}(\left\langle a \right\rangle, C_1)[\bar{A}_{k_1+1}] \simeq \text{Cay}(\left\langle a \right\rangle, C_2)[\bar{A}_{l_1+1}]$.
    \par Conversely, by Proposition \ref{pro: weakly connected component to another Cay}, $\text{Cay}(\left\langle a \right\rangle, C_1)[\bar{A}_{k_1}]$, $\text{Cay}(\left\langle a \right\rangle, C_2)[\bar{A}_{l_1}]$ can be regarded as Cayley graphs with connection set $C'_1 = \{b^{\frac{i_1}{d}}, b^{\frac{i_2}{d}}, \cdots, b^{\frac{i_n}{d}}\}$ and $C'_2 = \{b^{\frac{j_1}{d}}, b^{\frac{j_2}{d}}, \cdots, b^{\frac{j_n}{d}}\}$ respectively. By \Cref{theorem: main}, we have that for all $k \in \left\{1,2,\cdots, n\right\}$
    \begin{align*}
        \frac{i_k}{i_1} \equiv \frac{j_k}{j_1} \pmod {\frac{r}{d}}.
    \end{align*}
    By Corollary \ref{coro: same connected components}, $\gcd(C_1,\frac{r}{d}) = \gcd(C_2, \frac{r}{d})$. Then $\gcd(C_1,r) = \gcd(C_2, r)$ and $Par(C_1, \rho_1,\rho_2) = Par(C_2, \rho_1, \rho_2)$. In $\text{Cay}(\left\langle a \right\rangle,C_1)$, the class $\bar{A}_{k_1}$ containing $a^{i_1}$ has one fewer element than the class $\bar{A}_{k_1+1}$ containing $a^{i_1+1}$. Similarly, In $\text{Cay}(\left\langle a \right\rangle,C_2)$, the class $\bar{A}_{l_1}$ containing $a^{j_1}$ has one fewer element than the class $\bar{A}_{l_1+1}$ containing $a^{j_1+1}$. For all $t \in \{0, 1, \cdots, d-1\}$, let $f_t$ be the labelled colour-permuting isomorphism from $\bar{A}_{k_1+t}$ to $\bar{A}_{l_1+t}$ and let $\varphi$ be the bijection from $C_1$ to $C_2$ given by $\varphi(a^{i_k}) = a^{j_k}$. 
    \par Define $f: \left\langle a \right\rangle \to \left\langle a \right\rangle$. $f(a^u) = f_t(a^u)$, if $a^u \in \bar{A}_{k_1+t}$. It is clear that $\text{Cay}(\left\langle a \right\rangle,C_1) \simeq \text{Cay}(\left\langle a \right\rangle,C_2)$. 
\end{proof}

\vspace{1 em}
In the above proof, we assumed that there are two distinct weakly connected components. If, up to isomorphism, there is only one type of weakly connected component, the theorem still holds.
\par In \Cref{the: many}, the condition additionally requires two connection sets to satisfy $\gcd(C_1, r) = \gcd(C_2,r) = \gcd(i_\alpha, r) = \gcd(j_\beta,r)$ for some $a^{i_\alpha} \in C_1$ and $a^{j_\beta} \in C_2$. Next, we assume $\gcd(C_1, r) \ne \gcd(i_\alpha, r) $ for any $a^{i_\alpha} \in C_1$ and $\gcd(C_2, r) \ne \gcd(i_\beta, r) $ for any $a^{j_\beta} \in C_2$.

\begin{theorem}\label{theorem: 3}
Let $C_1 = \left\{a^{i_1}, a^{i_2}, \cdots, a^{i_n}\right\} \subseteq \textbf{CP}, C_2=\left\{a^{j_1}, a^{j_2}, \cdots, a^{j_n}\right\} \subseteq \textbf{CP}$ and $\gcd(C_1, r) \ne \gcd(i_\alpha, r) $ for any $a^{i_\alpha} \in C_1$ and $\gcd(C_2, r) \ne \gcd(i_\beta, r) $ for any $a^{j_\beta} \in C_2$. Then $\text{Cay}(\left\langle a \right\rangle,C_1) \simeq \text{Cay}(\left\langle a \right\rangle,C_2)$ if and only if there exist $a^{i_t} \in \textbf{CP} \backslash C_1$ and $a^{j_t} \in \textbf{CP} \backslash C_2$ satisfy $\gcd(i_t,i_1, i_2, \cdots, i_n, r) = \gcd(i_t, r) = \gcd(j_t,j_1, j_2, \cdots, j_n, r) = \gcd(j_t, r)$, such that $\text{Cay}(\left\langle a \right\rangle,C_1\cup \left\{a^{i_t}\right\}) \simeq \text{Cay}(\left\langle a \right\rangle, \allowbreak C_2\cup \left\{a^{j_t}\right\})$. 
\end{theorem}
\begin{proof}
    Suppose that $a^{i_t} \in \textbf{\textit{CP}} \backslash C_1$ and $a^{j_t} \in \textbf{\textit{CP}} \backslash C_2$ satisfy $\gcd(i_t,i_1, i_2, \cdots, \allowbreak i_n, r) = \gcd(i_t, r) = \gcd(j_t,j_1, j_2, \cdots, j_n, r) = \gcd(j_t, r)$, such that $\text{Cay}(\left\langle a \right\rangle, \allowbreak C_1\cup \left\{a^{i_t}\right\}) \simeq  \text{Cay}(\left\langle a \right\rangle,C_2\cup \left\{a^{j_t}\right\})$. Since the arcs generated by $a^{i_t}$ only be mapped to the arcs generated by $a^{j_t}$, $\text{Cay}(\left\langle a \right\rangle,C_1) \simeq  \text{Cay}(\left\langle a \right\rangle,C_2)$. 
    \par Conversely, suppose that $\text{Cay}(\left\langle a \right\rangle,C_1) \simeq  \text{Cay}(\left\langle a \right\rangle,C_2)$. Then there exist bijections $f : \left\langle a \right\rangle \to \left\langle a \right\rangle$ and $\varphi : C_1 \to C_2$ such that for all $a \in \left\langle a \right\rangle \text{, } x \in C_1$, $f(xa) = \varphi(x)f(a)$. Let $d = \gcd(i_1, i_2, \cdots, i_n, r)$, then by the $\textit{B\'ezout } Lemma$, there exist $x_1,x_2, \cdots, x_n, \bar{x} \in \mathbb{Z}$ such that $d = x_1i_1 +x_2i_2 + \cdots+ x_ni_n+\bar{x}r$. That means $a^d = a^{x_1i_1 +x_2i_2 +\cdots+ x_ni_n+\bar{x}r} = a^{x_1i_1 +x_2i_2 +\cdots+ x_ni_n}$. For arbitrary $x_t$, $x_t$ may be a  negative, suppose that there $\alpha_t>0 $, such that $y_t = x_t + \alpha_t r >0$. Let $M =\alpha_1i_1 +\alpha_2i_2+ \cdots +\alpha_ni_n$, then
    \begin{gather*}
        a^{d + Mr} = a^{x_1i_1 +x_2i_2 +\cdots+ x_ni_n} = \\
        a^{(x_1+\alpha_1r)i_1 +(x_2+\alpha_2r)i_2 +\cdots+ (x_n+\alpha_nr)i_n} = a^{y_1i_1+ y_2i_2+ \cdots +y_ni_n}.
    \end{gather*}
    That means $a^{d + Mr} = a^{y_1i_1+ y_2i_2+ \cdots +y_ni_n} \in \left\langle C_1 \right\rangle$. Let $a^{d+Mr} = a^{i_t}$ and $a^{i_t} = a^{i_{k_1}} \cdot a^{i_{k_2}} \cdot \cdots \cdot a^{i_{k_\alpha}}$, $a^{j_t} = \varphi(a^{i_{k_1}}) \cdot \varphi(a^{i_{k_2}}) \cdot \cdots \cdot \varphi(a^{i_{k_\alpha}})$. For arbitrary $x \in \left\langle a \right\rangle$ and $a^s \in C_1\cup \left\{a^{i_t}\right\}$, if $a^s \in C_1$, $f(a^sx) = \varphi(a^s)f(x)$; if $a^s = a^{i_t}$, $f(a^{i_t}x) = f(a^{i_{k_1}} \cdot a^{i_{k_2}} \cdot \cdots \cdot a^{i_{k_\alpha}}x) = \varphi(a^{i_{k_1}}) \cdot \varphi(a^{i_{k_2}}) \cdot \cdots \cdot \varphi(a^{i_{k_\alpha}}) f(x) = a^{j_t}f(x)$. Let 
    \begin{align*}
        \bar{\varphi}(a^x) = 
        \begin{cases}
            a^{j_t}, & \text{if }a^x = a^{i_t}, \\
            \varphi(a^x), & \text{if }a^x \in C_1.
        \end{cases}
    \end{align*}
    It follows that $\bar{\varphi}$ is also a bijection from $C_1\cup \left\{a^{i_t}\right\}$ to $C_2\cup \left\{a^{j_t}\right\}$. In this way, we can obtain $f(a^{i_t}x) = \bar{\varphi}(a^{i_t})f(x)$. Thus 
    \begin{align*}
        \text{Cay}(\left\langle a \right\rangle,C_1\cup \left\{a^{i_t}\right\}) \simeq  \text{Cay}(\left\langle a \right\rangle,C_2\cup \left\{a^{j_t}\right\}).
    \end{align*}
    \par Since $i_t \equiv Mr+d \pmod r $, $i_t \equiv d \pmod r$ can be obtained. Let $d = i_t +\beta r $, $d|i_t$ can be obtained by $d|r$. Let $r = dm$, then $i_t = d + \beta m =(1+\beta m)d$. Thus $d$ is a common factor of $ i_1, i_2, \cdots, i_n, i_t, r$. Suppose that $c$ is also a common factor of $ i_1, i_2, \cdots, i_n, i_t, r$, then $c|d$ can be obtained by $c$ is a common factor of $ i_1, i_2, \cdots, i_n, r$. Thus $\gcd(i_1, i_2, \cdots, i_n, i_t, r) = d$. It is clear that $\gcd(1+\beta m, m) = 1$, then $\gcd(d(1+\beta m), dm)= \gcd(i_t, r)=d$. Thus $\gcd(i_1, i_2, \cdots, i_n, i_t, r) = \gcd(i_t, r)$. 
    \par Since $\gcd(j_t,r)=d<\gcd(j_\beta,r)$ for every $a^{j_\beta}\in C_2$, $a^{j_t}$ is the only element of $C_2\cup \left\{a^{j_t}\right\}$ whose greatest common divisor with $r$ equals $d$; the cycle-length argument therefore forces $\varphi(a^{i_t})=a^{j_t}$. Thus $\text{Cay}(\left\langle a \right\rangle,C_1) \simeq  \text{Cay}(\left\langle a \right\rangle,C_2) \Longleftrightarrow \text{Cay}(\left\langle a \right\rangle,C_1\cup \left\{a^{i_t}\right\}) \simeq \text{Cay}(\left\langle a \right\rangle,C_2\cup \left\{a^{j_t}\right\})$.
\end{proof}

Likewise, under this circumstance, the simple case of a connection set of size $1$ can also be derived.

\begin{corollary}\label{coro: single2}
Let $C_1 = \left\{a^{i}\right\} \subseteq \textbf{\textit{CP}}, C_2=\left\{a^{j}\right\} \subseteq \textbf{\textit{CP}}$. $\text{Cay}(\left\langle a \right\rangle,C_1) \simeq  \text{Cay}(\left\langle a \right\rangle,C_2)$ if and only if 
$C_1$ and $C_2$ satisfy the following conditions,
\begin{flalign*}
    &(1) \text{ } \gcd(i, r) = \gcd(j ,r).&\\
    &(2) \text{ There exists a rotation } \phi, \text{ such that }
    \phi(X_1) = Y_1.
\end{flalign*}
\end{corollary}

\begin{proposition}\label{pro: cp ncp}
Suppose that $C_1 = \left\{a^{i_1}, a^{i_2}, \cdots, a^{i_n}\right\} \subseteq \left\langle a \right\rangle$, $i_1,i_2,...,i_k< m-1; i_{k+1}, i_{k+2}, \cdots ,i_n \ge m-1$ and $C_2 =\left\{a^{j_1}, a^{j_2}, \cdots, a^{j_k}, a^{j_{k+1}}, \cdots, a^{j_n}\right\} \subseteq \left\langle a \right\rangle$, $j_1 = i_1,j_2=i_2,...,j_k = i_k< m-1,j_{k+1}, j_{k+2}, \cdots ,j_n \ge m-1$. Let $\text{Cay}(\left\langle a \right\rangle,C_1) \simeq \text{Cay}(\left\langle a \right\rangle,C_2)$ and $f$ be labelled colour-permuting isomorphism, then $f$ is the graph automorphism on $\text{Cay}(\left\langle a \right\rangle, \left\{a^{i_t}\right\})$ for all $i_t, j_t < m-1$.
\end{proposition}
\begin{proof}
    Suppose that $f$ is the labelled colour-permuting isomorphism and $\varphi$ is the bijection from $C_1$ to $C_2$. Take $a^{i_t} \in C_1$, and $1 \le t \le k$, for all $a^t \in \left\langle a \right\rangle$,  $f(a^{i_t}a^t) = \varphi(a^{i_t}) f(a^t) \Longleftrightarrow (f(a^t), f(a^{i_t}a^t)) = (f(a^t), \varphi(a^{i_t}) f(a^t)) \Longleftrightarrow \text{Cay}(\left\langle a \right\rangle, \left\{a^{i_t}\right\}) \simeq \text{Cay}(\left\langle a \right\rangle, \left\{a^{j_t}\right\})$ and $f$ is the labelled colour-permuting isomorphism. By Corollary \ref{coro: 2}, $\varphi(a^{i_t}) = a^{j_t} \in \textbf{\textit{NCP}}$ for all $i_t < m-1$ and $\varphi(a^{i_t}) = a^{j_t} = a^{i_t}$. Thus $f$ is the graph automorphism on $\text{Cay}(\left\langle a \right\rangle, \left\{a^{i_t}\right\})$ for all $i_t, j_t < m-1$.
\end{proof}

\begin{corollary}\label{coro: 1}
Let $C_1 = \left\{a^{i_1}, a^{i_2}, \cdots, a^{i_n}\right\} \subseteq \left\langle a \right\rangle$, $C_2 =\left\{a^{j_1}, a^{j_2}, \cdots, a^{j_n}\right\} \subseteq \left\langle a \right\rangle$ and $C_1\cap \textbf{CP} \ne \emptyset, C_1\cap \textbf{NCP} \ne \emptyset, C_2\cap \textbf{CP} \ne \emptyset, C_2\cap \textbf{NCP} \ne \emptyset$. If $\text{Cay}(\left\langle a \right\rangle,C_1) \simeq \text{Cay}(\left\langle a \right\rangle,C_2)$ and $f$ is a labelled colour-permuting isomorphism, then $f$ is also a labelled colour-permuting isomorphism from $\text{Cay}(\left\langle a \right\rangle, C_1 \cap \textbf{CP})$ to $\text{Cay}(\left\langle a \right\rangle, C_2 \cap \textbf{CP})$.
\end{corollary}

\begin{theorem}\label{theorem: 4} 
Suppose that $C_1 = \left\{a^{i_1}, a^{i_2}, \cdots, a^{i_n}\right\} \subseteq \left\langle a \right\rangle$, $i_1,i_2,...,i_k< m-1,i_{k+1}, \cdots ,i_n \ge m-1$, $C_2 =\left\{a^{j_1}, a^{j_2}, \cdots, a^{j_k}, a^{j_{k+1}}, \cdots, a^{j_n}\right\} \subseteq \left\langle a \right\rangle$, $j_1 = i_1,j_2=i_2,...,j_k = i_k< m-1,j_{k+1}, j_{k+2}, \cdots ,j_n \ge m-1$. Then $\text{Cay}(\left\langle a \right\rangle,C_1) \simeq \text{Cay}(\left\langle a \right\rangle,C_2)$ if and only if $\text{Cay}(\left\langle a \right\rangle, C_1 \cap \textbf{CP}) \simeq \text{Cay}(\left\langle a \right\rangle, \allowbreak C_2 \cap \textbf{CP})$, and there exists a labelled colour-permuting isomorphism $f$ from $\text{Cay}(\left\langle a \right\rangle, C_1 \cap \textbf{CP})$ to $\text{Cay}(\left\langle a \right\rangle, C_2 \cap \textbf{CP})$ such that it is a graph automorphism on $\text{Cay}(\left\langle a \right\rangle, \left\{a^{i_t}\right\})$ for all $i_t< m-1$.
\end{theorem}
\begin{proof}
    Suppose that $\text{Cay}(\left\langle a \right\rangle,C_1) \simeq \text{Cay}(\left\langle a \right\rangle,C_2)$. By Proposition~\ref{pro: cp ncp}, $f$ is a labelled colour-permuting isomorphism from $\text{Cay}(\left\langle a \right\rangle, \left\{a^{i_t}\right\})$ to $\text{Cay}(\left\langle a \right\rangle, \left\{a^{j_t}\right\})$ for all $i_t, j_t < m-1$. And by the corollary~\ref{coro: 1}, $f$ is also a labelled colour-permuting isomorphism from $\text{Cay}(\left\langle a \right\rangle, C_1 \cap \textbf{\textit{CP}})$ to $\text{Cay}(\left\langle a \right\rangle, C_2 \cap \textbf{\textit{CP}})$.
    \par Conversely, it is clear that $f$ is also the labelled colour-permuting isomorphism from $\text{Cay}(\left\langle a \right\rangle,C_1)$ to $\text{Cay}(\left\langle a \right\rangle,C_2)$.
\end{proof}

The index-period decomposition yields an effective framework for studying Cayley graphs on finite monogenic semigroups. In this paper, we described the weakly connected components in terms of the arithmetic data of the connection set and gave a complete characterization of labelled colour-permuting isomorphisms in the non-cyclic, cyclic, and mixed cases. The resulting criteria show that the cyclic part governs the rotational behavior of the graph, while the non-cyclic part imposes rigid restrictions on the possible isomorphisms. These results indicate that the CI problem for finite monogenic semigroups differs essentially from the classical group case and point to further study of CI monogenic semigroups and related semigroup classes.

\end{document}